# STABILITY RESULTS FOR SYSTEMS DESCRIBED BY COUPLED RETARDED FUNCTIONAL DIFFERENTIAL EQUATIONS AND FUNCTIONAL DIFFERENCE EQUATIONS

**Iasson Karafyllis[*], Pierdomenico Pepe[**] and Zhong-Ping Jiang[***]**

[*] Department of Environmental Engineering, Technical University of Crete,
73100, Chania, Greece
email: ikarafyl@enveng.tuc.gr

[**]Dipartimento di Ingegneria Elettrica, Universita degli Studi dell'Aquila,
Monteluco di Roio, 67040, L'Aquila, Italy
email: pepe@ing.univaq.it

[***]Department of Electrical and Computer Engineering, Polytechnic University,
Six Metrotech Center, Brooklyn, NY 11201, U.S.A.
email: zjiang@control.poly.edu

**Abstract**

In this work stability results for systems described by coupled Retarded Functional Differential Equations (RFDEs) and Functional Difference Equations (FDEs) are presented. The results are based on the observation that the composite system can be regarded as the feedback interconnection of a subsystem described by RFDEs and a subsystem described by FDEs. Recent Small-Gain results and Lyapunov-like characterizations of the Weighted Input-to-Output Stability property for systems described by RFDEs and FDEs are employed. It is shown that the stability results provided in this work can be used to study stability for systems described by neutral functional differential equations and systems described by hyperbolic partial differential equations.

## 1. Introduction

In this work we consider control systems described by coupled Retarded Functional Differential Equations (RFDEs) and Functional Difference Equations (FDEs). Let $D \subseteq \Re^l$ be a non-empty set, $U \subseteq \Re^m$ be a non-empty set with $0 \in U$ and consider the system described by the following equations:

$$\dot{x}_1(t) = f_1(t, d(t), T_{r_1}(t)x_1, T_{r_2-\tau(t)}(t-\tau(t))x_2, u(t)) \tag{1.1a}$$

$$x_2(t) = f_2(t, d(t), T_{r_1}(t)x_1, T_{r_2-\tau(t)}(t-\tau(t))x_2, u(t)) \tag{1.1b}$$

$$Y(t) = H(t, T_{r_1}(t)x_1, T_{r_2}(t)x_2, u(t)) \in Y$$
$$x_1(t) \in \Re^{n_1}, x_2(t) \in \Re^{n_2}, d(t) \in D, u(t) \in U, t \geq 0 \tag{1.1c}$$

where $r_1, r_2 \geq 0$, $f_i : \bigcup_{t \geq 0} \{t\} \times D \times C^0([-r_1,0];\Re^{n_1}) \times L^\infty([-r_2+\tau(t),0];\Re^{n_2}) \times U \to \Re^{n_i}$, $i=1,2$, $H : \Re^+ \times C^0([-r_1,0];\Re^{n_1}) \times L^\infty([-r_2,0];\Re^{n_2}) \times U \to Y$ ($Y$ is a normed linear space) are locally bounded mappings with $f_i(t,d,0,0,0) = 0$ $i=1,2$, $H(t,0,0,0) = 0$ for all $(t,d) \in \Re^+ \times D$. Specifically, we consider systems of the form (1.1) with initial conditions $x_1(t_0+\theta) = x_{10}(\theta); \theta \in [-r_1, 0]$ and $x_2(t_0+\theta) = x_{20}(\theta); \theta \in [-r_2, 0]$ with $x_{10} \in C^0([-r_1,0]; \Re^{n_1})$, $x_{20} \in L^\infty([-r_2,0]; \Re^{n_2})$, under the following hypotheses:

**(P1)** The function $\tau: \Re^+ \to (0, +\infty)$ is continuous with $\sup_{t \geq 0} \tau(t) \leq r_2$.

**(P2)** There exist functions $a \in K_\infty$, $\beta \in K^+$ such that
$$\left|f_i(t,d,x_1, T_{r_2-\tau(t)}(-\tau(t))x_2, u)\right| \leq a\left(\beta(t)\|x_1\|_{r_1}\right) + a\left(\beta(t)\|T_{r_2-\tau(t)}(-\tau(t))x_2\|_{r_2-\tau(t)}\right) + a(\beta(t)|u|) \quad i = 1, 2, \text{ for all}$$
$(t, d, x_1, x_2, u) \in \Re^+ \times D \times C^0([-r_1, 0]; \Re^{n_1}) \times L^\infty([-r_2, 0]; \Re^{n_2}) \times U$.

**(P3)** For every $x_1 \in C^0([-r_1, +\infty); \Re^{n_1})$, $d \in L^\infty_{loc}(\Re^+; D)$, $u \in L^\infty_{loc}(\Re^+; U)$ and $x_2 \in L^\infty_{loc}([-r_2, +\infty); \Re^{n_2})$ the mappings $t \to f_i(t, d(t), T_{r_1}(t)x_1, T_{r_2-\tau(t)}(t-\tau(t))x_2, u(t))$, $i = 1, 2$ are measurable. Moreover, for each fixed $(t, d, x_2, u) \in \Re^+ \times D \times L^\infty([-r_2, 0]; \Re^{n_2}) \times U$ the mapping $f_1(t, d, x_1, T_{r_2-\tau(t)}(-\tau(t))x_2, u)$ is continuous with respect to $x_1 \in C^0([-r_1, 0]; \Re^{n_1})$.

**(P4)** For every pair of bounded sets $I \subset \Re^+$ and $\Omega \subset C^0([-r_1, 0]; \Re^{n_1}) \times L^\infty([-r_2, 0]; \Re^{n_2}) \times U$, there exists $L := L(I, \Omega) \geq 0$ such that

$$(x_1(0) - y_1(0))'\left(f_1(t, d, x_1, T_{r_2-\tau(t)}(-\tau(t))x_2, u) - f_1(t, d, y_1, T_{r_2-\tau(t)}(-\tau(t))x_2, u)\right) \leq L\|x_1 - y_1\|_{r_1}^2 \quad (1.2)$$
$$\forall (t, d) \in I \times D, \forall (x_1, x_2, u) \in \Omega, \forall (y_1, x_2, u) \in \Omega$$

**(P5)** The mapping $H: \Re^+ \times C^0([-r_1, 0]; \Re^{n_1}) \times L^\infty([-r_2, 0]; \Re^{n_2}) \times U \to Y$ is continuous with $H(t, 0, 0, 0) = 0$ for all $t \geq 0$. Moreover, the image set $H(\Omega)$ is bounded for each bounded set $\Omega \subset \Re^+ \times C^0([-r_1, 0]; \Re^{n_1}) \times L^\infty([-r_2, 0]; \Re^{n_2}) \times U$.

For example hypotheses (P1), (P2), (P3) are satisfied if $D \subset \Re^l$ is compact and there exist continuous functions $\tau_i: \Re^+ \to (0, +\infty)$ ($i = 1, ..., p$), $\tau: \Re^+ \to (0, +\infty)$ with $0 < \tau_1(t) < \tau_2(t) < ... < \tau_p(t) \leq \tau(t)$ for all $t \geq 0$ and $\sup_{t \geq 0} \tau(t) \leq r_2$, continuous mappings $g_i: \Re^+ \times D \times C^0([-r_1, 0]; \Re^{n_1}) \times \Re^{pn_2} \times \Re^k \times U \to \Re^{n_i}$, $i = 1, 2$, $h: \Re^+ \times [-r_2, 0] \times \Re^{n_2} \to \Re^k$ with $g_i(t, d, 0, 0, 0, 0) = 0$, $h(t, \theta, 0) = 0$ for all $(t, \theta, d) \in \Re^+ \times [-r-T, 0] \times D$, such that

$$f_i(t, d, x_1, x_2, u) = g_i\left(t, d, x_1, x_2(-\tau_1), x(-\tau_2), ..., x_2(-\tau_p), \int_{-r_2}^{-\tau(t)} h(t, \theta, x_2(\theta))d\theta, u\right), \quad i = 1, 2$$
for all $(t, d, x_1, x_2, u) \in \Re^+ \times D \times C^0([-r_1, 0]; \Re^{n_1}) \times L^\infty([-r_2, 0]; \Re^{n_2}) \times U$.

The reason for allowing the output to take values in abstract normed linear spaces is that the case (1.1) allows the study of:

- outputs with no delays, e.g. $Y(t) = h(t, x_1(t), x_2(t))$ with $Y = \Re^k$,
- outputs with discrete or distributed delay, e.g. $Y(t) = h(x_1(t), x_1(t-r_1), x_2(t), x_2(t-r_2))$ or $Y(t) = \int_{t-r_1}^{t} h(t, \theta, x_1(\theta))d\theta$ with $Y = \Re^k$,
- functional outputs with memory, e.g. $Y(t) = h(t, \theta, x_1(t+\theta)); \theta \in [-r_1, 0]$ or the identity output $Y(t) = \begin{cases} x_1(t+\theta); \theta \in [-r_1, 0] \\ x_2(t+\theta); \theta \in [-r_2, 0] \end{cases}$ with $Y = C^0([-r_1, 0]; \Re^{n_1}) \times L^\infty([-r_2, 0]; \Re^{n_2})$.

**Motivation for the study of system (1.1):**

Systems of the form (1.1) arise in many problems in Mathematical Control Theory and Mathematical Systems Theory (see for instance [4,5] and the references therein). For example, consider the stabilization problem for the scalar system:

$$\dot{x}(t) = f(x(t)) + u(t) + a\,u(t-r)$$
$$x(t) \in \Re, u(t) \in \Re \tag{1.3}$$

where $f: \Re \to \Re$ is a continuous function with $f(0) = 0$ and $r > 0$, $a \in \Re$ are constants. Notice the way that the input $u(t)$ appears in the equation (1.3). If the designer selects to apply the feedback linearization approach for system (1.3), then we have:

$$u(t) = -Kx(t) - f(x(t)) - au(t-r) \tag{1.4}$$

where $K > 0$. Consequently, if $a \neq 0$ the closed-loop system (1.3) with (1.4) is described by the following system of coupled RFDEs and FDEs:

$$\dot{x}(t) = -Kx(t)$$
$$u(t) = -Kx(t) - f(x(t)) - au(t-r) \tag{1.5}$$
$$x(t) \in \Re, u(t) \in \Re$$

Notice that system (1.5) has the form of system (1.1) with $u(t)$ in place of $x_2(t)$ and $x(t)$ in place of $x_1(t)$. Moreover, hypotheses (P1-4) are satisfied for system (1.5).

However, the strongest motive for the study of systems of the form (1.1) is that systems of the form (1.1) allow the consideration of **discontinuous solutions** to systems described by Neutral Functional Differential Equations. For example, consider the scalar system described by a Neutral Functional Differential Equation:

$$\frac{d}{dt}(x(t) - x(t-2r)) = x(t-r), \; x(t) \in \Re \tag{1.6}$$

with initial condition $T_{2r}(t_0) = x_0 \in C^0([-2r, 0]; \Re)$. The solution of (1.6) for $t \in [t_0, t_0 + r]$ is given by:

$$x(t) = x(t-2r) + x(t_0) - x(t_0 - 2r) + \int_{t_0 - r}^{t-r} x(s)\,ds \tag{1.7}$$

It is clear from (1.7) that the solution of (1.6) can be defined even if the initial condition is discontinuous, i.e., $T_{2r}(t_0) = x_0 \in L^\infty([-2r, 0]; \Re)$. How to obtain such a (weak) solution?

The following idea was proposed in [2] for linear systems (though it was expressed in a different way): First define

$$x_1(t) = x(t) - x(t-2r)$$
$$x_2(t) = x(t)$$

Then equation (1.6) is transformed to the following system of coupled RFDEs and FDEs of the form (1.1):

$$\dot{x}_1(t) = x_2(t-r)$$
$$x_2(t) = x_1(t) + x_2(t-2r) \tag{1.8}$$
$$x_1(t) \in \Re, x_2(t) \in \Re$$

Notice that hypotheses (P1-4) are satisfied for system (1.8) (hypothesis (P5) is irrelevant since there is no output). The solution of (1.8) for $t \in (t_0, t_0 + r]$ is given by:

$$x_1(t) = x_1(t_0) + \int_{t_0-r}^{t-r} x_2(s)ds$$

$$x_2(t) = x_2(t-2r) + x_1(t_0) + \int_{t_0-r}^{t-r} x_2(s)ds \tag{1.9}$$

Notice that $x_2(t) = x(t)$ does not coincide with the solution of (1.6) given by (1.7) unless $x_1(t_0) = x_2(t_0) - x_2(t_0 - 2r)$, the so-called "matching condition". It should be emphasized that if the "matching condition" does not hold then the solution of (1.8) (given by (1.9)) is discontinuous, even if the initial condition is smooth. Consequently, system (1.8) provides a generalized framework for the study of the Neutral Functional Differential Equation (1.6).

The idea described for the simple example (1.6) can be generalized for nonlinear control systems described by Neutral Functional Differential Equations of the following form (special case of the so-called Hale's form, see [6]):

$$\frac{d}{dt}\left(x(t) - g(t, T_{r-\tau(t)}(t-\tau(t))x)\right) = f(t, d(t), T_r(t)x, u(t)) \,,\, x(t) \in \Re^n \tag{1.10}$$

Without loss of generality we may assume that the continuous function $\tau: \Re^+ \to (0,+\infty)$ with $\sup_{t\geq 0} \tau(t) \leq r$, is non-increasing. If we define $x_1(t) = x(t) - g(t, T_{r-\tau(t)}(t-\tau(t))x)$, $x_2(t) = x(t)$ and the operator:

$$\Re^+ \times C^0\left([-r,0]; \Re^{n_1}\right) \times L^\infty\left([-2r,0]; \Re^{n_2}\right) \ni (t, x_1, x_2) \to G(t, x_1, T_{2r-\tau(t)}(-\tau(t))x_2)$$

$$G(t, x_1, T_{2r-\tau(t)}(-\tau(t))x_2) := \begin{cases} x_2(\theta), & \theta \in [-r, -\tau(t)] \\ x_1(\theta) + g(t+\theta, T_{r-\tau(t+\theta)}(\theta - \tau(t+\theta))x_2), & \theta \in (-\tau(t), 0] \end{cases}$$

then system (1.10) is associated with the following system described by coupled RFDEs and FDEs:

$$\dot{x}_1(t) = f\left(t, d(t), G\left(t, T_r(t)x_1, T_{2r-\tau(t)}(t-\tau(t))x_2\right), u(t)\right)$$

$$x_2(t) = x_1(t) + g(t, T_{r-\tau(t)}(t-\tau(t))x_2) \tag{1.11}$$

Notice that system (1.11) is a system of the form (1.1). The component $x_2$ of the solution of (1.11) coincides with the solution $x$ of (1.10) if and only if the initial data are continuous functions which satisfy the "matching condition": $G\left(t_0, T_r(t_0)x_1, T_{2r-\tau(t_0)}(t_0 - \tau(t_0))x_2\right) = T_r(t_0)x_2$. However, notice that even if the matching condition holds the solution of (1.11) can be defined for discontinuous initial data. Consequently, if the matching condition holds and the initial data are discontinuous then the component $x_2$ of the solution of (1.11) is a discontinuous mapping which satisfies the differential equation $\frac{d}{dt}\left(x_2(t) - g(t, T_{r-\tau(t)}(t-\tau(t))x_2)\right) = f(t, d(t), T_r(t)x_2, u(t))$, almost everywhere for $t \geq t_0$. Thus, if the matching condition holds, then system (1.11) provides "weak" solutions to the Neutral Functional Differential Equation (1.10). Other concepts of weak solutions for linear neutral functional differential equations were given in [2,9]. In recent works control-theoretic aspects for linear neutral functional differential equations are studied (see [7,8]).

The approach described above is not restricted to Neutral Functional Differential Equations of the form (1.10). We can also consider Neutral Functional Differential Equations of the form (Bellman's form, see [1]):

$$\dot{x}(t) = f(t, d(t), T_r(t)x, T_{r-\tau(t)}(t-\tau(t))\dot{x}, u(t)) \,,\, x(t) \in \Re^n \tag{1.12}$$

In this case the corresponding system of coupled RFDEs and FDEs is:

$$\dot{x}_1(t) = f\left(t, d(t), T_r(t)x_1, T_{r-\tau(t)}(t-\tau(t))x_2, u(t)\right)$$
$$\dot{x}_2(t) = f\left(t, d(t), T_r(t)x_1, T_{r-\tau(t)}(t-\tau(t))x_2, u(t)\right) \quad (1.13)$$
$$x_1(t) \in \Re^n, x_2(t) \in \Re^n, d(t) \in D, u(t) \in U, t \geq 0$$

Notice that if $T_r(t_0)x_2 = T_r(t_0)\dot{x}$ and $T_r(t_0)x_1 = T_r(t_0)x$ then the component $x_1$ of the solution of (1.13) coincides with the solution $x$ of (1.12) for all $t \geq t_0$.

Consequently, it should be emphasized that the study of coupled RFDEs and FDEs offers a great advantage: Neutral Functional Differential Equations of the form (1.10) and Neutral Functional Differential Equations of the form (1.12) can be studied **in the same way and in the same framework**.

Another field which motivates the study of systems of the form (1.1) is the field of control (or dynamical) systems described by hyperbolic partial differential equations of the form:

$$\frac{\partial v_i}{\partial t}(t,z) + a_i \frac{\partial v_i}{\partial z}(t,z) = f_i(t, v_i(t,z), \xi(t)), \quad i = 1,...,p \quad (1.14a)$$
$$v_i(t,z) \in \Re^{n_i}, \xi(t) \in \Re^k, z \in [0,1]$$

where $a_i > 0$ ($i = 1,...,p$) are constants, along with boundary conditions of the form:

$$v_i(t,0) = F_i\left(t, d(t), \xi(t), u(t), v(t,z); z \in [\tau(t),1]\right), \quad i = 1,...,p \quad (1.14b)$$

where $v(t,z) = (v_1(t,z),...,v_p(t,z))'$ and

$$\dot{\xi}(t) = g\left(t, d(t), \xi(t), u(t), v(t,z); z \in [\tau(t),1]\right) \quad (1.14c)$$

The problem (1.14) is accompanied with initial conditions $v_i(t_0,z) = v_{0i}(z)$ and $\xi(t_0) = \xi_0 \in \Re^k$. Initial value problems of the form (1.14) arise in electrical, thermal and hydraulic engineering (see for instance the model of combined heat and electricity generation in [29] and other models reported in [21] concerning lossless transmission lines with electrical circuits and turbines under waterhammer conditions).

If we define $x_1(t) = \xi(t)$ and $x_2(t) = (v_1(t,0),...,v_p(t,0))'$, then it can be shown that the state variables $x_1(t)$, $x_2(t)$ satisfy a system of coupled RFDEs and FDEs for $t \geq t_0 + \max_{i=1,...,p} a_i^{-1}$. Consequently, the asymptotic behaviour of system (1.14) is determined by the associated system of coupled RFDEs and FDEs. The discontinuous solutions generated by the associated system of coupled RFDEs and FDEs are important, since such solutions correspond to "weak" solutions of the problem (1.14).

Finally, it should be noticed that recent contributions in the literature study systems of coupled RFDEs and FDEs of the form (1.1) per se (see [4,19,21,23,24,25,26,28,29]).

In this work we present sufficient Lyapunov-like conditions for (Uniform) Robust Global Asymptotic Output Stability (RGAOS) and (Uniform) Weighted Input-to-Output Stability (WIOS) for systems of the form (1.1). Our results are based **on the decomposition of system (1.1) as the feedback interconnection of a system described by RFDEs and a system described by FDEs**. This particular viewpoint allows us to study the stability properties of (1.1) in great generality as well as to obtain a unified framework for a wide class of stability notions, including the notion of Input-to-State Stability (ISS). It should be emphasized that the introduction of the notion of ISS by E. D. Sontag in [30,31,32] for finite-dimensional systems described by ordinary differential equations, led to an exceptionally rich period of progress in mathematical systems and control theory. The notion of ISS was extended to the notion of Input-to-Output Stability (IOS) in [33,34,35] and to non-uniform in time notions of ISS and IOS in [13,14,15,16] (which extended the applicability of ISS to time-varying systems). It is our belief that the notions of ISS and IOS have become one of the most important conceptual tools for the development of nonlinear robust stability and control theory for a wide class of dynamical systems and consequently, **one of the novelties of the present work is the study of ISS/IOS for system (1.1)**.

The present work is structured as follows: in Section 2, we provide some preliminary results on existence and uniqueness of solutions of (1.1), which allow us to consider system (1.1) as a control system with a robust equilibrium point, in the sense described in [15]. In Section 3, we present a stability result, which is based on the Small-Gain Theorem given in [17]. The result relies on the notion of (Uniform) Weight Input-to-Output Stability for control systems with outputs. Sufficient Lyapunov-like and Razumikhin-like conditions are also presented. In Section 4, examples are presented where the stability results of Section 3 are utilized. Finally, in Section 5 we provide the concluding remarks of this work.

**Notations** Throughout this paper we adopt the following notations:

* For a vector $x \in \Re^n$ we denote by $|x|$ its usual Euclidean norm and by $x'$ its transpose. For a bounded function $x:[-r,0] \to \Re^n$ we define $\|x\|_r := \sup_{\theta \in [-r,0]} |x(\theta)|$. For a matrix $A \in \Re^{m \times n}$ by $|A|$ we denote the induced norm of the matrix, i.e., $|A| := \sup\{|Ax| ; x \in \Re^n, |x| = 1\}$. $I \in \Re^{n \times n}$ denotes the identity matrix.

* Let $I \subseteq \Re$ be an interval. By $C^0(I;\Omega)$, we denote the class of continuous functions on $I$, which take values in $\Omega \subseteq \Re^n$. By $C^1(I;\Omega)$, we denote the class of functions on $I$ with continuous derivative, which take values in $\Omega \subseteq \Re^n$. By $L^\infty(I;\Omega)$ ($L^\infty_{loc}(I;\Omega)$), we denote the class of measurable and (locally) bounded functions on $I$, which take values in $\Omega \subseteq \Re^n$. If $\Omega \subseteq \Re^n$ is a subspace of $\Re^n$, $L^\infty(I;\Omega)$ ($L^\infty_{loc}(I;\Omega)$) is a normed linear space with norm $\sup_{t \in I} |x(t)|$, for $x \in L^\infty(I;\Omega)$ ($x \in L^\infty_{loc}(I;\Omega)$).

* $\Re^+$ denotes the set of non-negative real numbers.

* We denote by $K^+$ the class of positive $C^0$ functions defined on $\Re^+$. We say that a function $\rho : \Re^+ \to \Re^+$ is positive definite if $\rho(0) = 0$ and $\rho(s) > 0$ for all $s > 0$. We say that a function $\rho : \Re^+ \to \Re^+$ is of class $N$, if $\rho$ is non-decreasing with $\rho(0) = 0$. By $K$ we denote the set of positive definite, increasing and continuous functions. We say that a positive definite, increasing and continuous function $\rho : \Re^+ \to \Re^+$ is of class $K_\infty$ if $\lim_{s \to +\infty} \rho(s) = +\infty$. By $KL$ we denote the set of all continuous functions $\sigma = \sigma(s,t) : \Re^+ \times \Re^+ \to \Re^+$ with the properties: (i) for each $t \geq 0$ the mapping $\sigma(\cdot,t)$ is of class $K$ ; (ii) for each $s \geq 0$, the mapping $\sigma(s,\cdot)$ is non-increasing with $\lim_{t \to +\infty} \sigma(s,t) = 0$.

* Let $x:[a-r,b) \to \Re^n$ with $b > a > -\infty$ and $r > 0$. By $T_r(t)x$ we denote the "$r$-history" of $x$ at time $t \in [a,b)$, i.e., $T_r(t)x := x(t+\theta) ; \theta \in [-r,0]$.

* By $\| \|_Y$, we denote the norm of the normed linear space $Y$.

## 2. Preliminary Results for Control Systems Described by Coupled RFDEs and FDEs

In this section we provide some fundamental results that allow us to consider system (1.1) under hypotheses (P1-5) as a control system with a robust equilibrium point in the sense described in [15]. We start with an existence-uniqueness-continuation theorem for the solution of (1.1). We say that a mapping $x:[a,b) \to \Re^n$ with $-\infty < a < b \leq +\infty$ is absolutely continuous on $[a,b)$ if for every $c \in (a,b)$ the mapping $x:[a,b) \to \Re^n$ is absolutely continuous on $[a,c]$.

**Theorem 2.1:** *Consider system (1.1) under hypotheses (P1-4). Then for every $t_0 \geq 0$, $(x_{10}, x_{20}) \in C^0([-r_1,0];\Re^{n_1}) \times L^\infty([-r_2,0];\Re^{n_2})$, $d \in L^\infty_{loc}(\Re^+;D)$, $u \in L^\infty_{loc}(\Re^+;U)$ there exists $t_{\max} \in (t_0,+\infty]$ and a unique pair of mappings $x_1 \in C^0([t_0-r_1,t_{\max});\Re^{n_1})$, $x_2 \in L^\infty_{loc}([t_0-r_2,t_{\max});\Re^{n_2})$ with $T_{r_1}(t_0)x_1 = x_{10}$, $T_{r_2}(t_0)x_2 = x_{20}$, $x_1 \in C^0([t_0-r_1,t_{\max});\Re^{n_1})$ being absolutely continuous on $[t_0,t_{\max})$ such that (1.1a) holds a.e. for $t \in [t_0,t_{\max})$ and (1.1b) holds for all $t \in (t_0,t_{\max})$. In addition, if $t_{\max} < +\infty$ then for every $M > 0$ there exists $t \in [t_0,t_{\max})$ with $\|T_{r_1}(t)x_1\|_{r_1} > M$.*

Theorem 2.1 guarantees that $t_{\max} \in (t_0, +\infty]$ is the maximal existence time for the solution of (1.1). The idea behind the proof of Theorem 2.1 is the method of steps, used already in [21,25].

**Proof of Theorem 2.1:** Let $t_0 \geq 0$, $(x_{10}, x_{20}) \in C^0\!\left([-r_1, 0]; \Re^{n_1}\right) \times L^\infty\!\left([-r_2, 0]; \Re^{n_2}\right)$, $d \in L^\infty_{loc}\!\left(\Re^+; D\right)$, $u \in L^\infty_{loc}\!\left(\Re^+; U\right)$ (arbitrary). Define $h := \min\!\left(1; \min\{\tau(t_0 + s) : s \in [0,1]\}\right)$. Notice that by virtue of definition of $h > 0$, it holds that $t - \tau(t) \leq t_0$ for all $t \in [t_0, t_0 + h]$. By virtue of Theorem 2.1 in [6] (and its extension for Caratheodory conditions in page 58 of the same book) there exists $\delta \in (0, h]$ and $x_1 \in C^0\!\left([t_0 - r_1, t_0 + \delta); \Re^{n_1}\right)$ with $T_{r_1}(t_0)x_1 = x_{10}$, being absolutely continuous on $[t_0, t_0 + \delta)$ such that the differential equation $\dot{x}_1(t) = f_1(t, d(t), T_{r_1}(t)x_1, T_{r_2-\tau(t)}(t - \tau(t))x_2, u(t))$ is satisfied a.e. for $t \in [t_0, t_0 + \delta)$. Moreover, hypothesis (P4) guarantees that the mapping $x_1 \in C^0\!\left([t_0 - r_1, t_0 + \delta); \Re^{n_1}\right)$ is unique (see [14,18]). There exist two cases for the mapping $x_1 \in C^0\!\left([t_0 - r_1, t_0 + \delta); \Re^{n_1}\right)$:

a) if $\delta < h$ then for every $M > 0$ there exists $t \in [t_0, t_0 + \delta)$ with $\|T_{r_1}(t)x_1\|_{r_1} > M$.

b) if $\|T_{r_1}(t)x_1\|_{r_1}$ is bounded for all $t \in [t_0, t_0 + \delta)$ then the mapping $x_1 \in C^0\!\left([t_0 - r_1, t_0 + \delta); \Re^{n_1}\right)$ can be extended continuously in a unique way on $[t_0 - r_1, t_0 + h]$.

Next we consider the FDE $x_2(t) = f_2(t, d(t), T_{r_1}(t)x_1, T_{r_2-\tau(t)}(t - \tau(t))x_2, u(t))$. By virtue of hypotheses (P2), (P3) there exists a unique mapping $x_2 \in L^\infty_{loc}\!\left([t_0 - r_2, t_0 + \delta); \Re^{n_2}\right)$ with $T_{r_2}(t_0)x_2 = x_{20}$ such that the FDE $x_2(t) = f_2(t, d(t), T_{r_1}(t)x_1, T_{r_2-\tau(t)}(t - \tau(t))x_2, u(t))$ is satisfied for all $t \in (t_0, t_0 + \delta)$. Moreover, if $\|T_{r_1}(t)x_1\|_{r_1}$ is bounded for all $t \in [t_0, t_0 + \delta)$ and the mapping $x_1 \in C^0\!\left([t_0 - r_1, t_0 + \delta); \Re^{n_1}\right)$ can be extended continuously in a unique way on $[t_0 - r_1, t_0 + h]$ then similarly $x_2 \in L^\infty\!\left([t_0 - r_2, t_0 + \delta); \Re^{n_2}\right)$ can be extended on $[t_0 - r_2, t_0 + h]$ (notice that hypothesis (P2) implies that $x_2$ is bounded as long as $x_1$ is bounded).

If case (a) holds then define $t_{\max} = t_0 + \delta$ and the proof is complete. If case (b) holds all arguments can be repeated with $t_0 + h$ in place of $t_0$ (next step). We continue the same procedure of construction of the solution step-by-step. The procedure may be stopped after some steps (if case (a) is encountered) or may be continued indefinitely (if case (a) is never encountered). In the latter case for each step $i$ we obtain a pair of mappings $x_1 \in C^0\!\left([t_0 - r_1, t_{i+1}]; \Re^{n_1}\right)$, $x_2 \in L^\infty\!\left([t_0 - r_2, t_{i+1}]; \Re^{n_2}\right)$ with $T_{r_1}(t_0)x_1 = x_{10}$, $T_{r_2}(t_0)x_2 = x_{20}$, $x_1 \in C^0\!\left([t_0 - r_1, t_{i+1}]; \Re^{n_1}\right)$ being absolutely continuous on $[t_0, t_{i+1}]$, such that (1.1a) holds a.e. for $t \in [t_0, t_{i+1}]$ and (1.1b) holds for all $t \in (t_0, t_{i+1}]$, where the sequence $\{t_i\}_{i=0}^\infty$ satisfies $t_{i+1} = t_i + \min\!\left(1; \min\{\tau(t_i + s) : s \in [0,1]\}\right)$ for all $i = 0,1,2,\ldots$. Notice that the sequence $\{t_i\}_{i=0}^\infty$ is increasing and consequently $\lim t_i = \sup t_i$. The assumption that $L = \lim t_i = \sup t_i < +\infty$ implies that $t_{i+1} \geq t_i + \mu$, where $\mu = \min\!\left(1; \min\{\tau(s) : s \in [0, L+1]\}\right)$ for all $i = 0,1,2,\ldots$, which gives the contradiction $t_i \geq t_0 + (i-1)\mu$ for all $i = 1,2,\ldots$. It follows that $\lim t_i = \sup t_i = +\infty$.

The proof is complete. ◁

**Remark 2.2:** According to Theorem 2.1 above, Definition 2.1 in [15] and Definition 2.4 in [15], system (1.1) under hypotheses (P1-5) is a control system $\Sigma := (\mathcal{X}, Y, M_U, M_D, \phi, \pi, H)$ with outputs that satisfies the BIC property with state space $\mathcal{X} = C^0\!\left([-r_1, 0]; \Re^{n_1}\right) \times L^\infty\!\left([-r_2, 0]; \Re^{n_2}\right)$, output space $Y$, set of allowable control inputs $M_U = L^\infty_{loc}\!\left(\Re^+; U\right)$, set of allowable disturbances $M_D = L^\infty_{loc}\!\left(\Re^+; D\right)$ and set of sampling times $\pi(t_0, x_0, u, d) = [t_0, t_{\max})$, where $t_{\max} > t_0$ is the maximal existence time of the solution. Moreover, if a finite escape time occurs then the component $x_1$ of the solution of (1.1) must be unbounded (but $x_2$ may or may not be unbounded).

The following theorem guarantees that $(0,0) \in C^0\!\left([-r_1, 0]; \Re^{n_1}\right) \times L^\infty\!\left([-r_2, 0]; \Re^{n_2}\right)$ is a robust equilibrium point from the input $u$ (in the sense of Definition 2.6 in [15]) for system (1.1) under hypotheses (P1-4).

**Theorem 2.3:** *Consider system (1.1) under hypotheses (P1-4). Then for every $\varepsilon > 0$, $T, h \in \Re^+$ there exists $\delta := \delta(\varepsilon, T, h) > 0$ such that for all $(t_0, x_{10}, x_{20}) \in [0,T] \times C^0([-r_1,0]; \Re^{n_1}) \times L^\infty([-r_2,0]; \Re^{n_2})$, $(u,d) \in L^\infty_{loc}(\Re^+; U) \times L^\infty_{loc}(\Re^+; D)$ with $\|x_{10}\|_{r_1} + \|x_{20}\|_{r_2} + \sup_{t \geq 0}|u(t)| < \delta$ there exists $t_{max} \in (t_0 + h, +\infty]$ and a unique pair of mappings $x_1 \in C^0([t_0 - r_1, t_{max}); \Re^{n_1})$, $x_2 \in L^\infty_{loc}([t_0 - r_2, t_{max}); \Re^{n_2})$ with $T_{r_1}(t_0)x_1 = x_{10}$, $T_{r_2}(t_0)x_2 = x_{20}$, $x_1 \in C^0([t_0 - r_1, t_{max}); \Re^{n_1})$ being absolutely continuous on $[t_0, t_{max})$, such that (1.1a) holds a.e. for $t \in [t_0, t_{max})$, (1.1b) holds for all $t \in (t_0, t_{max})$ and*

$$\sup\left\{ \|T_{r_1}(t)x_1\|_{r_1} + \|T_{r_2}(t)x_2\|_{r_2} \; ; \; t \in [t_0, t_0 + h] \right\} < \varepsilon \tag{2.1}$$

**Proof:** The proof of Theorem 2.3 relies on the following fact, which is proved at the Appendix.

**Fact I:** *For every $\varepsilon > 0$, $T \in \Re^+$ there exists $\tilde{\delta} := \tilde{\delta}(\varepsilon, T) > 0$ such that for all $(t_0, x_{10}, x_{20}) \in [0,T] \times C^0([-r_1,0]; \Re^{n_1}) \times L^\infty([-r_2,0]; \Re^{n_2})$, $(u,d) \in L^\infty_{loc}(\Re^+; U) \times L^\infty_{loc}(\Re^+; D)$ with $\|x_{10}\|_{r_1} + \|x_{20}\|_{r_2} + \sup_{t \geq 0}|u(t)| < \tilde{\delta}$ there exists $t_{max} \in (t_0 + h, +\infty]$ where $h := h(T) = \min(1; \min\{\tau(s) : s \in [0, T+1]\})$ and a unique pair of mappings $x_1 \in C^0([t_0 - r_1, t_{max}); \Re^{n_1})$, $x_2 \in L^\infty_{loc}([t_0 - r_2, t_{max}); \Re^{n_2})$ with $T_{r_1}(t_0)x_1 = x_{10}$, $T_{r_2}(t_0)x_2 = x_{20}$ such that (1.1a) holds a.e. for $t \in [t_0, t_{max})$, (1.1b) holds for all $t \in (t_0, t_{max})$ and (2.1) holds.*

Next consider the sequence $\{T_i\}_{i=0}^\infty$ which is generated by the recursive relation:

$$T_{i+1} = T_i + \min(1; \min\{\tau(s) : s \in [0, T_i + 1]\}), \; i = 0,1,2,\ldots \text{ with } T_0 = T \geq 0 \tag{2.2}$$

As in the proof of Theorem 2.1 it can be shown (by contradiction) that $\lim T_i = \sup T_i = +\infty$ for all $T_0 \geq 0$. Consequently, given arbitrary $T, h \in \Re^+$, there exists non-negative integer $l(T,h)$ such that the sequence $\{T_i\}_{i=0}^\infty$ defined by (2.2) with initial condition $T_0 = T$, satisfies $T_i \geq T + h$ for all $i \geq l(T, h)$. The following fact exploits the properties of the sequence $\{T_i\}_{i=0}^\infty$ defined by (2.2).

**Fact II:** *For every $\varepsilon > 0$, $T \in \Re^+$, $i$ non-negative integer, there exists $\delta_i := \delta_i(\varepsilon, T) > 0$ such that for all $(t_0, x_{10}, x_{20}) \in [0,T] \times C^0([-r_1,0]; \Re^{n_1}) \times L^\infty([-r_2,0]; \Re^{n_2})$, $(u,d) \in L^\infty_{loc}(\Re^+; U) \times L^\infty_{loc}(\Re^+; D)$ with $\|x_{10}\|_{r_1} + \|x_{20}\|_{r_2} + \sup_{t \geq 0}|u(t)| < \delta_i$ there exists $t_{max} \in (t_0 + h, +\infty]$, where $h := h(T) = T_{i+1} - T$, $\{T_i\}_{i=0}^\infty$ is the sequence that satisfies (2.2), and a unique pair of mappings $x_1 \in C^0([t_0 - r_1, t_{max}); \Re^{n_1})$, $x_2 \in L^\infty_{loc}([t_0 - r_2, t_{max}); \Re^{n_2})$ with $T_{r_1}(t_0)x_1 = x_{10}$, $T_{r_2}(t_0)x_2 = x_{20}$ such that (1.1a) holds a.e. for $t \in [t_0, t_{max})$, (1.1b) holds for all $t \in (t_0, t_{max})$ and (2.1) holds.*

The proof of Fact II will be made by induction. By virtue of Fact I it is clear that Fact II holds for $i = 0$. Suppose that Fact II holds for certain non-negative integer $i$. Let $\varepsilon > 0$, $T \in \Re^+$ and define:

$$\delta_{i+1}(\varepsilon, T) := \min\left\{ \delta_i(\varepsilon, T) \; ; \; \delta_i\left(\frac{1}{2}\tilde{\delta}(\varepsilon, T_{i+1}), T\right); \frac{1}{2}\tilde{\delta}(\varepsilon, T_{i+1}) \right\} > 0 \tag{2.3}$$

Next consider the solution $x_1 \in C^0([t_0 - r_1, t_{max}); \Re^{n_1})$, $x_2 \in L^\infty_{loc}([t_0 - r_2, t_{max}); \Re^{n_2})$ of (1.1) with $T_{r_1}(t_0)x_1 = x_{10}$, $T_{r_2}(t_0)x_2 = x_{20}$ corresponding to inputs $(u,d) \in L^\infty_{loc}(\Re^+; U) \times L^\infty_{loc}(\Re^+; D)$ with $\|x_{10}\|_{r_1} + \|x_{20}\|_{r_2} + \sup_{t \geq 0}|u(t)| < \delta_{i+1}$. Since $\delta_{i+1} \leq \delta_i$, it follows from the assumption that Fact II holds for the non-negative integer $i$:

$$\sup\left\{ \|T_{r_1}(t)x_1\|_{r_1} + \|T_{r_2}(t)x_2\|_{r_2} \; ; \; t \in [t_0, t_0 + T_{i+1} - T] \right\} < \varepsilon \tag{2.4}$$

Moreover, since $\delta_{i+1}(\varepsilon,T) \leq \delta_i\left(\frac{1}{2}\tilde{\delta}(\varepsilon,T_{i+1}),T\right)$, it follows from the assumption that Fact II holds for the non-negative integer $i$:

$$\left\|T_{r_1}(t_0+T_{i+1}-T)x_1\right\|_{r_1} + \left\|T_{r_2}(t_0+T_{i+1}-T)x_2\right\|_{r_2} < \frac{1}{2}\tilde{\delta}(\varepsilon,T_{i+1}) \tag{2.5}$$

Furthermore, since $\|x_{10}\|_{r_1} + \|x_{20}\|_{r_2} + \sup_{t\geq 0}|u(t)| < \delta_{i+1}$ and $\delta_{i+1}(\varepsilon,T) \leq \frac{1}{2}\tilde{\delta}(\varepsilon,T_{i+1})$, we obtain that $\sup_{t\geq 0}|u(t)| \leq \frac{1}{2}\tilde{\delta}(\varepsilon,T_{i+1})$. Combining (2.5) and the previous inequality we get:

$$\left\|T_{r_1}(t_0+T_{i+1}-T)x_1\right\|_{r_1} + \left\|T_{r_2}(t_0+T_{i+1}-T)x_2\right\|_{r_2} + \sup_{t\geq 0}|u(t)| < \tilde{\delta}(\varepsilon,T_{i+1}) \tag{2.6}$$

Notice that since $t_0 \in [0,T]$, we obtain that $t_0+T_{i+1}-T \in [0,T_{i+1}]$. The solution $x_1 \in C^0\left([t_0-r_1,t_{\max});\Re^{n_1}\right)$, $x_2 \in L^\infty_{loc}\left([t_0-r_2,t_{\max});\Re^{n_2}\right)$ of (1.1) with $T_{r_1}(t_0)x_1 = x_{10}$, $T_{r_2}(t_0)x_2 = x_{20}$ corresponding to inputs $(u,d) \in L^\infty_{loc}(\Re^+;U) \times L^\infty_{loc}(\Re^+;D)$ with $\|x_{10}\|_{r_1} + \|x_{20}\|_{r_2} + \sup_{t\geq 0}|u(t)| < \delta_{i+1}$ coincides for $t \geq t_0+T_{i+1}-T$ with the solution $x_1 \in C^0\left([t_0+T_{i+1}-T-r_1,t_{\max});\Re^{n_1}\right)$, $x_2 \in L^\infty_{loc}\left([t_0+T_{i+1}-T-r_2,t_{\max});\Re^{n_2}\right)$ of (1.1) with initial condition $(T_{r_1}(t_0+T_{i+1}-T)x_1, T_{r_2}(t_0+T_{i+1}-T)x_2)$ corresponding to same inputs $(u,d) \in L^\infty_{loc}(\Re^+;U) \times L^\infty_{loc}(\Re^+;D)$ satisfying (2.6). Using Fact I, in conjunction with (2.6), definition (2.2) and the fact that $t_0+T_{i+1}-T \in [0,T_{i+1}]$, we obtain:

$$\sup\left\{\left\|T_{r_1}(t)x_1\right\|_{r_1} + \left\|T_{r_2}(t)x_2\right\|_{r_2} ; \; t \in [t_0+T_{i+1}-T, t_0+T_{i+2}-T]\right\} < \varepsilon \tag{2.7}$$

Combining (2.4) with (2.7), we may conclude that Claim II holds for $i+1$.

By virtue of Fact II, it follows that Theorem 2.3 holds with $\delta(\varepsilon,T,h) := \delta_{l(T,h)}(\varepsilon,T) > 0$, where $l(T,h)$ is the non-negative integer with the property that the sequence $\{T_i\}_{i=0}^\infty$ defined by (2.2) with initial condition $T_0 = T$, satisfies $T_i \geq T + h$ for all $i \geq l(T,h)$. The proof is complete. ◁

**Remark 2.4:** It should be emphasized that Theorems 2.1 and 2.3 guarantee that all stability results obtained in [13,15] for general control systems with the "Boundedness Implies Continuation" property (BIC property, see [15]) hold as well for system (1.1) under hypotheses (P1-5). This implication enables us to obtain the stability results of the following section.

**Remark 2.5:** It is important to notice that Theorems 2.1 and 2.3 can be applied to systems described by FDEs of the form:

$$\begin{aligned} &x(t) = f(t,d(t),T_{r-\tau(t)}(t-\tau(t))x,u(t)) \\ &x(t) \in \Re^n, d(t) \in D, u(t) \in U, t \geq 0 \end{aligned} \tag{2.8}$$

where $D \subseteq \Re^l$ is a non-empty set, $U \subseteq \Re^m$ is a non-empty set with $0 \in U$, $r > 0$, $f: \bigcup_{t\geq 0}\{t\} \times D \times L^\infty([-r+\tau(t),0];\Re^n) \times U \to \Re^n$, under the following hypotheses:

**(Q1)** The function $\tau: \Re^+ \to (0,+\infty)$ is continuous with $\sup_{t\geq 0}\tau(t) \leq r$.

**(Q2)** There exist functions $a \in K_\infty$, $\beta \in K^+$ such that $|f(t,d,T_{r-\tau(t)}(-\tau(t))x,u)| \leq a\left(\beta(t)\|T_{r-\tau(t)}(-\tau(t))x\|_{r-\tau(t)}\right) + a(\beta(t)|u|)$, for all $(t,d,x,u) \in \Re^+ \times D \times L^\infty([-r,0];\Re^n) \times U$.

**(Q3)** For every $d \in L^\infty_{loc}(\Re^+;D)$, $u \in L^\infty_{loc}(\Re^+;U)$ and $x \in L^\infty_{loc}([-r,+\infty);\Re^n)$ the mapping $t \to f(t,d(t),T_{r-\tau(t)}(t-\tau(t))x,u(t))$ is measurable.

Indeed, system (2.8) can be embedded into the following system described by coupled RFDEs and FDEs:

$$\dot{\xi}(t) = 0$$
$$x(t) = f(t, d(t), T_{r-\tau(t)}(t-\tau(t))x, u(t)) \quad (2.9)$$
$$\xi(t) \in \Re, x(t) \in \Re^n, d(t) \in D, u(t) \in U, t \geq 0$$

which is a system of the form (1.1) that satisfies hypotheses (P1-4). Consequently, Theorems 2.1 and 2.3 can be applied to system (2.9) and we obtain:

**Corollary 2.6:** *Consider system (2.8) under hypotheses (Q1-3). Then for every $t_0 \geq 0$, $x_0 \in L^\infty([-r,0];\Re^n)$, $(u,d) \in L^\infty_{loc}(\Re^+;U) \times L^\infty_{loc}(\Re^+;D)$ there exists a unique mapping $x \in L^\infty_{loc}([t_0-r,+\infty);\Re^n)$ with $T_r(t_0)x = x_0$, such that (2.8) holds for all $t > t_0$. Moreover, for every $\varepsilon > 0$, $T,h \in \Re^+$ there exists $\delta := \delta(\varepsilon,T,h) > 0$ such that for all $(t_0,x_0) \in [0,T] \times L^\infty([-r,0];\Re^n)$, $(u,d) \in L^\infty_{loc}(\Re^+;U) \times L^\infty_{loc}(\Re^+;D)$ with $\|x_0\|_r + \sup_{t \geq 0}|u(t)| < \delta$ the solution $x(t)$ of (2.8) with initial condition $T_r(t_0)x = x_0$, corresponding to inputs $(u,d) \in L^\infty_{loc}(\Re^+;U) \times L^\infty_{loc}(\Re^+;D)$ satisfies $\sup\{\|T_r(t)x\|_r\;;\;t \in [t_0,t_0+h]\;\} < \varepsilon$.*

## 3. Stability Results

In this section we present stability results for a wide class of systems described by coupled RFDEs and FDEs. Particularly, we consider the following class of systems described by coupled RFDEs and FDEs:

$$\dot{x}_1(t) = f_1(t, d(t), T_{r_1}(t)x_1, u(t), H_2(t, T_{r_2-\tau(t)}(t-\tau(t))x_2))) \quad (3.1a)$$

$$x_2(t) = f_2(t, d(t), T_{r_2-\tau(t)}(t-\tau(t))x_2, u(t), H_1(t, T_{r_1}(t)x_1))$$
$$x_1(t) \in \Re^{n_1}, x_2(t) \in \Re^{n_2}, u(t) \in U, d(t) \in D, t \geq 0 \quad (3.1b)$$

$$Y(t) = H(t, T_{r_1}(t)x_1, T_{r_2}(t)x_2) \in Y \quad (3.1c)$$

where $r_1, r_2 \geq 0$, $D \subseteq \Re^l$ a non-empty set, $U \subseteq \Re^m$ a non-empty set with $0 \in U$, $Y$ is a normed linear space, $H_1 : \Re^+ \times C^0([-r_1,0];\Re^{n_1}) \to S_1$, $H_2 : \cup_{t \geq 0}\{t\} \times L^\infty([-r_2+\tau(t),0];\Re^{n_2}) \to S_2$, $H : \Re^+ \times C^0([-r_1,0];\Re^{n_1}) \times L^\infty([-r_2,0];\Re^{n_2}) \to Y$ are continuous mappings, $S_1 \subseteq \Re^{k_1}$, $S_2 \subseteq \Re^{k_2}$ are sets with $0 \in S_1$, $0 \in S_2$ and the mappings $f_1 : \Re^+ \times D \times C^0([-r_1,0];\Re^{n_1}) \times U \times S_2 \to \Re^{n_1}$, $f_2 : \cup_{t \geq 0}\{t\} \times D \times L^\infty_{loc}([-r_2-\tau(t),0];\Re^{n_2}) \times U \times S_1 \to \Re^{n_2}$ are locally bounded mappings, which satisfy the following hypotheses:

**(R1)** The function $\tau : \Re^+ \to (0,+\infty)$ is continuous with $\sup_{t \geq 0}\tau(t) \leq r_2$.

**(R2)** For every $v \in L^\infty_{loc}(\Re^+;S_1)$, $d \in L^\infty_{loc}(\Re^+;D)$, $u \in L^\infty_{loc}(\Re^+;U)$ and $x_2 \in L^\infty_{loc}([-r_2,+\infty);\Re^{n_2})$ the mapping $t \to f_2(t, d(t), T_{r_2-\tau(t)}(t-\tau(t))x_2, u(t), v(t))$ is measurable.

**(R3)** The output map $H_1 : \Re^+ \times C^0([-r_1,0];\Re^{n_1}) \to S_1$, is a continuous mapping that maps bounded sets of $\Re^+ \times C^0([-r_1,0];\Re^{n_1})$ into bounded sets of $\Re^{k_1}$ with $H_1(t,0) = 0$ for all $t \geq 0$.

**(R4)** The mapping $H_2 : \cup_{t \geq 0}\{t\} \times L^\infty([-r_2+\tau(t),0];\Re^{n_2}) \to S_2$ is a continuous mapping that maps bounded subsets of $\cup_{t \geq 0}\{t\} \times L^\infty([-r_2+\tau(t),0];\Re^{n_2})$ into bounded sets of $\Re^{k_2}$ with $H_2(t,0) = 0$ for all $t \geq 0$. Moreover, for every $x_2 \in L^\infty_{loc}([-r_2,+\infty);\Re^{n_2})$ the mapping $t \to H_2(t, T_{r_2-\tau(t)}(t-\tau(t))x_2)$ is measurable.

**(R5)** There exist functions $\beta \in K^+$, $a \in K_\infty$ such that
$$\left| f_2\left(t, d, T_{r_2-\tau(t)}(-\tau(t))x_2, u, v\right) \right| \leq a\left( \beta(t) \sup_{\theta \in [-r_2, -\tau(t)]} |x_2(\theta)| \right) + a(\beta(t)|u|) + a(\beta(t)|v|)$$
for all $(t, x_2, u, v, d) \in \Re^+ \times L^\infty([-r_2, 0]; \Re^{n_2}) \times U \times S_1 \times D$ and $|f_1(t, d, x, u, v)| \leq a(\beta(t)\|x\|_{r_1}) + a(\beta(t)|u|) + a(\beta(t)|v|)$ for all $(t, x, u, v, d) \in \Re^+ \times C^0([-r_1, 0]; \Re^n) \times U \times S_2 \times D$.

**(R6)** The mapping $(x, u, v, d) \to f_1(t, d, x, u, v)$ is continuous for each fixed $t \geq 0$ and such that for every bounded $I \subseteq \Re^+$ and for every bounded $\Omega \subset C^0([-r_1, 0]; \Re^n) \times U \times S_2$, there exists a constant $L \geq 0$ such that:

$$(x(0) - y(0))'\left(f_1(t, d, x, u, v) - f_1(t, d, y, u, v)\right) \leq L \max_{\tau \in [-r_1, 0]} |x(\tau) - y(\tau)|^2$$

$\forall t \in I$, $\forall (x, u, v, y, u, v) \in \Omega \times \Omega$, $\forall d \in D$

**(R7)** There exists a countable set $A \subset \Re^+$, which is either finite or $A = \{t_k \ ; k = 1, ..., \infty\}$ with $t_{k+1} > t_k > 0$ for all $k = 1, 2, ...$ and $\lim t_k = +\infty$, such that mapping $(t, x, u, v, d) \in (\Re^+ \setminus A) \times C^0([-r, 0]; \Re^n) \times U \times S_2 \times D \to f_1(t, d, x, u, v)$ is continuous. Moreover, for each fixed $(t_0, x, u, v, d) \in \Re^+ \times C^0([-r, 0]; \Re^n) \times U \times S_2 \times D$, we have $\lim_{t \to t_0^+} f_1(t, d, x, u, v) = f_1(t_0, d, x, u, v)$.

**(R8)** The mapping $H : \Re^+ \times C^0([-r_1, 0]; \Re^{n_1}) \times L^\infty([-r_2, 0]; \Re^{n_2}) \to Y$ is continuous with $H(t, 0, 0) = 0$ for all $t \geq 0$. Moreover, the image set $H(\Omega)$ is bounded for each bounded set $\Omega \subset \Re^+ \times C^0([-r_1, 0]; \Re^{n_1}) \times L^\infty([-r_2, 0]; \Re^{n_2})$.

By virtue of Lemma 3.2 in [13] and Lemma 1 (page 4) in [3], it follows that system (3.1) under hypotheses (R1-8) is a system of the form (1.1) which satisfies hypotheses (P1-5). However, it should be emphasized that not every system of the form (1.1) can be expressed in the form (3.1).

Next, we consider the following systems:

$$\begin{aligned}
\dot{x}_1(t) &= f_1(t, d(t), T_{r_1}(t)x_1, u(t), v_1(t)) \\
Y_1(t) &= H_1(t, T_{r_1}(t)x_1) \\
x_1(t) &\in \Re^{n_1}, \ Y_1(t) \in S_1, \ (u(t), v_1(t)) \in U \times S_2, \ d(t) \in D, \ t \geq 0
\end{aligned} \quad (3.2)$$

which is a system described by RFDEs and the following system described by FDEs:

$$\begin{aligned}
x_2(t) &= f_2(t, d(t), T_{r_2-\tau(t)}(t-\tau(t))x_2, u(t), v_2(t)) \\
Y_2(t) &= H_2(t, T_{r_2-\tau(t)}(t-\tau(t))x_2) \\
x_2(t) &\in \Re^{n_2}, \ Y_2(t) \in S_2, \ (u(t), v_2(t)) \in U \times S_1, \ d(t) \in D, \ t \geq 0
\end{aligned} \quad (3.3)$$

The following things can be remarked for systems (3.2), (3.3):

\* The theory of retarded functional differential equations guarantees that under hypotheses (R3-7), for each $(t_0, x_{10}) \in \Re^+ \times C^0([-r_1, 0]; \Re^{n_1})$ and for each triple of measurable and locally bounded inputs $v_1 \in L^\infty_{loc}(\Re^+; S_2)$, $d \in L^\infty_{loc}(\Re^+; D)$, $u \in L^\infty_{loc}(\Re^+; U)$ there exists a unique absolutely continuous mapping $x_1(t)$ that satisfies a.e. the differential equation (3.2) with initial condition $T_{r_1}(t_0)x_1 = x_{10} \in C^0([-r_1, 0]; \Re^{n_1})$ (see [18]). Moreover, Theorem 3.2 in [6] (page 46) guarantees that (3.2) is a control system $\Sigma_1 := (C^0([-r_1, 0]; \Re^{n_1}), \Re^{k_1}, M_{U \times S_2}, M_D, \phi, \pi, H_1)$ with outputs that satisfies the Boundedness Implies Continuation property with $M_{U \times S_2}, M_D$ the sets of all measurable and locally bounded mappings $(u, v) : \Re^+ \to U \times S_2$, $d : \Re^+ \to D$, respectively (in the sense described in [15]). Furthermore, the classical semigroup property is satisfied for this system, i.e., we have

$\pi(t_0, x_0, u, d) = [t_0, t_{\max})$, where $t_{\max} > t_0$ is the maximal existence time of the solution. Finally, hypotheses (R3-7) guarantee that $0 \in C^0([-r_1, 0]; \Re^{n_1})$ is a robust equilibrium point from the input $(u, v) \in M_{U \times S_2}$ for $\Sigma_1$.

∗ Hypotheses (R1-5) guarantee that for each $(t_0, x_{20}) \in \Re^+ \times X$ with $X := L^\infty([-r_2, 0]; \Re^{n_2})$ and for each triple $v_2 \in L^\infty_{loc}(\Re^+; S_1)$, $d \in L^\infty_{loc}(\Re^+; D)$, $u \in L^\infty_{loc}(\Re^+; U)$ there exists a unique measurable and locally bounded mapping $x_2(t)$ that satisfies the difference equations (3.3) for all $t > t_0$ with initial condition $x_2(t_0 + \theta) = x_{20}(\theta); \theta \in [-r_2, 0]$. Consequently, (3.3) describes a control system $\Sigma_2 := (X, \Re^{k_2}, M_{U \times S_1}, M_D, \phi, \pi, H_2)$ with outputs, $M_{U \times S_1}$ being the set of all measurable and locally bounded functions $(u, v): \Re^+ \to U \times S_1$ and $M_D$ being the set of all measurable and locally bounded functions $d: \Re^+ \to D$ (in the sense described in [15]). Furthermore, Remark 2.5 and Corollary 2.6 show that system (3.3) is Robustly Forward Complete from the input $(u, v) \in M_{U \times S_1}$ and that $0 \in X$ is a robust equilibrium point from the input $(u, v) \in M_{U \times S_1}$ for system (3.3) in the sense described in [15]. Finally, notice that the classical semigroup property is satisfied for system (3.3), i.e., we have $\pi(t_0, x_{20}, u, d) = [t_0, +\infty)$.

Systems (3.1), (3.2) and (3.3) are related in the following way: it can be said that system (3.1) is the feedback interconnection of subsystems (3.2) and (3.3), in the sense described in [15]. Picture 1 presents schematically the interconnection of subsystems (3.2) and (3.3) that produces the composite system (3.1).

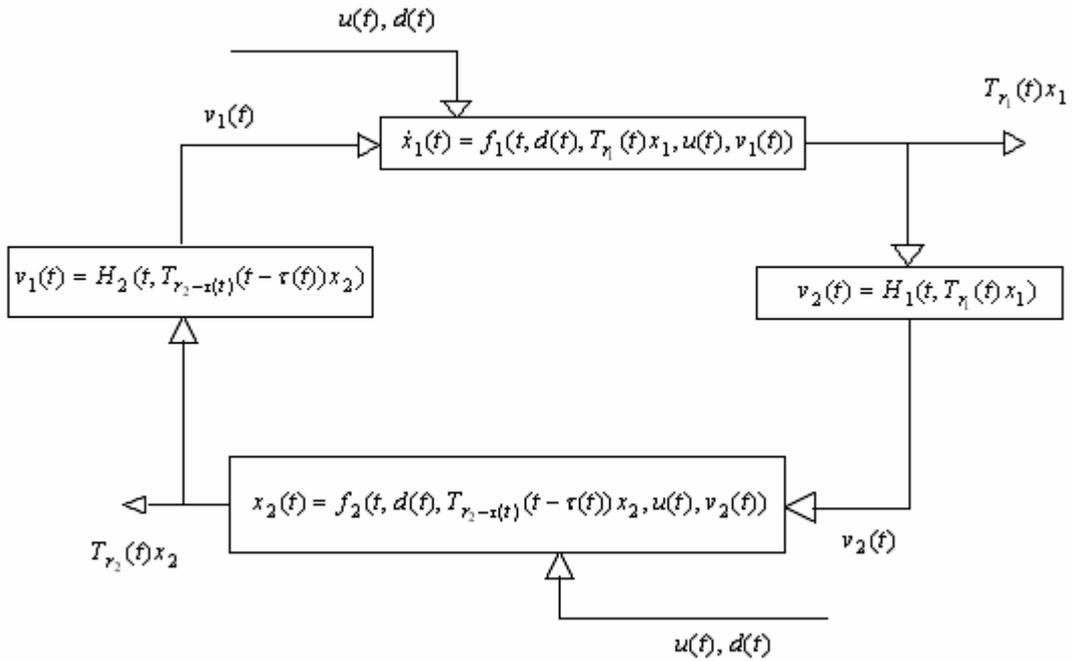

**Figure 1:** System (3.1) regarded as the feedback interconnection of subsystem (3.2) described by RFDEs and subsystem (3.3) described by FDEs

We are now in a position to present our main result, which is a direct of consequence of the Small-Gain Theorem presented in [17]. The definition of the (Uniform) Weighted Input-to-Output Stability ((U)WIOS) property is given in the Appendix for reader's convenience (see also [17]).

**Theorem 3.1:** *Consider system (3.1) under hypotheses (R1-8) and assume that:*

**(H1)** *Subsystem (3.2) satisfies the WIOS property from the inputs $v_1$ and $u$. Particularly, there exist functions $\sigma_1 \in KL$, $\beta_1, \mu_1, c_1, \delta_1, \delta_1^u, q_1^u \in K^+$, $\gamma_1, \gamma_1^u, a_1, p_1, p_1^u \in N$, such that the following estimate holds for all $(t_0, x_{10}, (v_1, u, d)) \in \Re^+ \times C^0([-r_1, 0]; \Re^{n_1}) \times L^\infty_{loc}(\Re^+; S_2 \times U \times D)$ and $t \geq t_0$ for the solution $x_1(t)$ of (3.2) with initial condition $T_{r_1}(t_0) x_1 = x_{10}$ corresponding to inputs $(v_1, u, d) \in L^\infty_{loc}(\Re^+; S_2 \times U \times D)$:*

$$\left| H_1(t, T_{r_1}(t) x_1) \right| \leq \sigma_1 \left( \beta_1(t_0) \| x_{10} \|_{r_1}, t - t_0 \right) + \sup_{t_0 \leq \tau \leq t} \gamma_1 \left( \delta_1(\tau) | v_1(\tau)| \right) + \sup_{t_0 \leq \tau \leq t} \gamma_1^u \left( \delta_1^u(\tau) | u(\tau)| \right) \tag{3.4}$$

$$\beta_1(t) \| T_{r_1}(t) x_1 \|_{r_1} \leq \max \left\{ \mu_1(t - t_0), c_1(t_0), a_1 \left( \| x_{10} \|_{r_1} \right), \sup_{t_0 \leq \tau \leq t} p_1 \left( | v_1(\tau)| \right), \sup_{t_0 \leq \tau \leq t} p_1^u \left( q_1^u(\tau) | u(\tau)| \right) \right\} \tag{3.5}$$

**(H2)** *Subsystem (3.3) satisfies the WIOS property from the inputs $v_2$ and $u$. Particularly, there exist functions $\sigma_2 \in KL$, $\beta_2, \mu_2, c_2, \delta_2, \delta_2^u, q_2^u \in K^+$, $\gamma_2, \gamma_2^u, a_2, p_2, p_2^u \in N$, such that the following estimate holds for all $(t_0, x_{20}, (v_2, u, d)) \in \Re^+ \times L^\infty([-r_2, 0]; \Re^{n_2}) \times L^\infty_{loc}(\Re^+; S_1 \times U \times D)$ and $t \geq t_0$ for the solution $x_2(t)$ of (3.3) with initial condition $T_{r_2}(t_0) x_2 = x_{20}$ corresponding to inputs $(v_2, u, d) \in L^\infty_{loc}(\Re^+; S_1 \times U \times D)$:*

$$\left| H_2(t, T_{r_2 - \tau(t)}(t - \tau(t)) x_2) \right| \leq \sigma_2 \left( \beta_2(t_0) \| x_{20} \|_{r_2}, t - t_0 \right) + \sup_{t_0 \leq s \leq t} \gamma_2 \left( \delta_2(s) | v_2(s)| \right) + \sup_{t_0 \leq s \leq t} \gamma_2^u \left( \delta_2^u(s) | u(s)| \right) \tag{3.6}$$

$$\beta_2(t) \| T_{r_2}(t) x_2 \|_{r_2} \leq \max \left\{ \mu_2(t - t_0), c_2(t_0), a_2 \left( \| x_{20} \|_{r_2} \right), \sup_{t_0 \leq s \leq t} p_2 \left( | v_2(s)| \right), \sup_{t_0 \leq s \leq t} p_2^u \left( q_2^u(s) | u(s)| \right) \right\} \tag{3.7}$$

**(H3)** *There exist functions $\rho \in K_\infty$, $a \in N$ and a constant $M > 0$ such that one of the following conditions holds for all $t, s \geq 0$ and $x = (x_1, x_2) \in C^0([-r_1, 0]; \Re^{n_1}) \times L^\infty([-r_2, 0]; \Re^{n_2})$:*

$$\delta_1(t) \leq M \tag{3.8a}$$

$$g_1 \left( \delta_1(t) g_2 \left( \max_{\theta \in [0, t]} \delta_2(\theta) s \right) \right) \leq s \tag{3.8b}$$

$$\| H(t, x_1, x_2) \|_Y \leq a \left( \left| H_1(t, x_1) \right| + \gamma_1 \left( \delta_1(t) \left| H_2(t, T_{r_2 - \tau(t)}(-\tau(t)) x_2) \right| \right) \right) \tag{3.8c}$$

*or*

$$\delta_2(t) \leq M \tag{3.9a}$$

$$g_2 \left( \delta_2(t) g_1 \left( \max_{\theta \in [0, t]} \delta_1(\theta) s \right) \right) \leq s \tag{3.9b}$$

$$\| H(t, x_1, x_2) \|_Y \leq a \left( \left| H_2(t, T_{r_2 - \tau(t)}(-\tau(t)) x_2) \right| + \gamma_2 \left( \delta_2(t) | H_1(t, x_1)| \right) \right) \tag{3.9c}$$

*where $g_i(s) := \gamma_i(s) + \rho(\gamma_i(s))$, $i = 1, 2$.*

*Then there exists a function $\gamma \in N$ such that system (3.1) satisfies the WIOS property from the input $u \in M_U$ with gain $\gamma \in N$ and weight $\delta \in K^+$, where*

$$\delta(t) := \max \{ \delta_1^u(t), \delta_2^u(t), q_1^u(t), q_2^u(t) \} \tag{3.10}$$

*Moreover, if $\beta_1, \beta_2, c_1, c_2, \delta_1, \delta_2 \in K^+$ are bounded then system (3.1) satisfies the UWIOS property from the input $u \in M_U$ with gain $\gamma \in N$ and weight $\delta \in K^+$.*

**Remark 3.2:**
  (a) It should be clear that Theorem 3.1 gives sufficient conditions (but not necessary) for the WIOS property for system (3.1). The main advantage of Theorem 3.1 is that the stability properties for system (3.1) can be verified by studying the stability properties of subsystems (3.2) and (3.3), which are simpler systems.
  (b) When $\gamma_1 \in N$ or $\gamma_2 \in N$ is identically zero, it follows that (3.8b) and (3.9b) are automatically satisfied. On the other hand, if $\gamma_i(s) = K_i s$ for certain constants $K_i \geq 0$ ($i=1,2$) then hypothesis (3.8b) (or (3.9b)) is satisfied if $K_1 K_2 \sup_{t \geq 0} \left( \delta_1(t) \max_{\tau \in [0,t]} \delta_2(\tau) \right) < 1$ (or $K_1 K_2 \sup_{t \geq 0} \left( \delta_2(t) \max_{\tau \in [0,t]} \delta_1(\tau) \right) < 1$).

In what follows, we present sufficient conditions so that subsystems (3.2) and (3.3) satisfy assumptions (H1) and (H2) of Theorem 3.1. The following theorem is a direct consequence of Theorem 4.6 in [18] and gives Lyapunov-like sufficient conditions that guarantee assumption (H1) for subsystem (3.1). Its proof can be found in the Appendix. We notice that for a functional $V : \Re^+ \times C^0([-r,0]; \Re^n) \to \Re$, the generalized derivative in the direction $v \in \Re^n$ is defined by

$$V^0(t,x;v) := \limsup_{\substack{h \to 0^+ \\ y \to 0, y \in C^0([-r,0];\Re^n)}} \frac{V(t+h, E_h(x;v) + hy) - V(t,x)}{h} \tag{3.11}$$

where $E_h(x;v)$ with $0 \leq h < r$, denotes the following operator:

$$E_h(x;v) := \begin{cases} x(0) + (\theta + h)v & \text{for } -h < \theta \leq 0 \\ x(\theta + h) & \text{for } -r \leq \theta \leq -h \end{cases} \tag{3.12}$$

Moreover, a continuous functional $V : \Re^+ \times C^0([-r,0]; \Re^n) \to \Re^+$ is "almost Lipschitz on bounded sets", if there exist non-decreasing functions $M : \Re^+ \to \Re^+$, $P : \Re^+ \to \Re^+$, $G : \Re^+ \to [1, +\infty)$ such that for all $R \geq 0$, the following properties hold:

**(1)** For every $x, y \in \{ x \in C^0([-r,0]; \Re^n) ; \|x\|_r \leq R \}$, it holds that: $|V(t,y) - V(t,x)| \leq M(R) \|y - x\|_r$, $\forall t \in [0, R]$

**(2)** For every absolutely continuous function $x : [-r,0] \to \Re^n$ with $\|x\|_r \leq R$ and essentially bounded derivative, it holds that:

$$|V(t+h,x) - V(t,x)| \leq hP(R)\left(1 + \sup_{-r \leq \tau \leq 0} |\dot{x}(\tau)|\right), \text{ for all } t \in [0, R] \text{ and } 0 \leq h \leq \frac{1}{G\left(R + \sup_{-r \leq \tau \leq 0} |\dot{x}(\tau)|\right)}$$

More details for "almost Lipschitz on bounded sets" functionals, as well as for the generalized derivative (3.11), are provided in [18].

**Theorem 3.3 (Lyapunov-like sufficient conditions for hypothesis (H1)):** *Consider system (3.2) under hypotheses (R3-7) and suppose that there exist a Lyapunov functional $V : \Re^+ \times C^0([-r_1, 0]; \Re^{n_1}) \to \Re^+$, which is almost Lipschitz on bounded sets, functions $a_2$ of class $K_\infty$, $\zeta, \zeta^u$ of class $N$, $\beta, \delta_1, \delta_1^u$ of class $K^+$ and a continuous positive definite function $\rho : \Re^+ \to \Re^+$ such that:*

$$V(t, x_1) \leq a_2\left(\beta(t) \|x_1\|_{r_1}\right), \quad \forall (t, x_1) \in \Re^+ \times C^0([-r_1, 0]; \Re^{n_1}) \tag{3.13}$$

$$V^0(t, x_1; f_1(t, d, x_1, u, v_1)) \leq -\rho(V(t, x_1)),$$
for all $(t, x_1, u, v_1, d) \in \Re^+ \times C^0([-r_1, 0]; \Re^{n_1}) \times U \times S_2 \times D$ with $\max\{\zeta(\delta_1(t)|v_1|); \zeta^u(\delta_1^u(t)|u|)\} \leq V(t, x_1)$ (3.14)

*Moreover, suppose that there exist functions $a_1, p$ of class $K_\infty$, $\mu$ of class $K^+$ and a constant $R \geq 0$ such that one of the following inequalities holds:*

$$a_1(|H_1(t,x_1)|) \leq V(t,x_1), \quad \forall (t,x_1) \in \Re^+ \times C^0([-r_1,0];\Re^{n_1}) \tag{3.15a}$$

or

$$p(\mu(t)|x_1(0)|) \leq V(t,x_1) + R, \quad \forall (t,x_1) \in \Re^+ \times C^0([-r_1,0];\Re^{n_1}) \tag{3.15b}$$

If

* (3.15a) holds and there exist functions $\mu_1, c_1, \phi \in K^+$, $g_1, p_1, p_1^u \in N$, such that for every $(t_0, x_{10}, (v_1, u, d)) \in \Re^+ \times C^0([-r_1,0];\Re^{n_1}) \times L_{loc}^\infty(\Re^+; S_2 \times U \times D)$ the solution $x_1(t)$ of (3.2) with initial condition $T_{r_1}(t_0)x_1 = x_{10}$ corresponding to inputs $(v_1, u, d) \in L_{loc}^\infty(\Re^+; S_2 \times U \times D)$ exists for all $t \geq t_0$ and satisfies the following estimate:

$$\phi(t)\|T_{r_1}(t)x_1\|_{r_1} \leq \max\left\{ \mu_1(t-t_0), c_1(t_0), g_1(\|x_{10}\|_{r_1}), \sup_{t_0 \leq \tau \leq t} p_1(|v_1(\tau)|), \sup_{t_0 \leq \tau \leq t} p_1^u(\delta_1^u(\tau)|u(\tau)|) \right\} \tag{3.16}$$

then there exists a function $\sigma_1 \in KL$, such that estimate (3.4) holds with $\beta_1(t) \equiv \beta(t)$, $\gamma_1(s) := a_1^{-1}(\zeta(s))$, $\gamma_1^u(s) := a_1^{-1}(\zeta^u(s))$ for all $(t_0, x_{10}, (v_1, u, d)) \in \Re^+ \times C^0([-r_1,0];\Re^{n_1}) \times L_{loc}^\infty(\Re^+; S_2 \times U \times D)$ and $t \geq t_0$ for the solution $x_1(t)$ of (3.2) with initial condition $T_{r_1}(t_0)x_1 = x_{10}$ corresponding to inputs $(v_1, u, d) \in L_{loc}^\infty(\Re^+; S_2 \times U \times D)$.

* (3.15b) holds and $\delta_1(t) \equiv 1$ then for every $(t_0, x_{10}, (v_1, u, d)) \in \Re^+ \times C^0([-r_1,0];\Re^{n_1}) \times L_{loc}^\infty(\Re^+; S_2 \times U \times D)$ the solution $x_1(t)$ of (3.2) with initial condition $T_{r_1}(t_0)x_1 = x_{10}$ corresponding to inputs $(v_1, u, d) \in L_{loc}^\infty(\Re^+; S_2 \times U \times D)$ exists for all $t \geq t_0$. Furthermore, for every $\phi \in K^+$ there exist functions $\mu_1, c_1 \in K^+$, $g_1, p_1, p_1^u \in N$, such that (3.16) holds for all $t \geq t_0$. Moreover, if $\phi \in K^+$ is bounded and there exists constant $L > 0$ such that:

$$\beta(t) + \frac{1}{\mu(t)} \leq L, \quad \forall t \geq 0 \tag{3.17}$$

then the function $c_1 \in K^+$ is bounded.

**Remark 3.4:** It should be emphasized that the main results in [18] show that the existence of a Lyapunov functional $V: \Re^+ \times C^0([-r_1,0];\Re^{n_1}) \to \Re^+$ satisfying the assumptions of Theorem 3.3 is a necessary and sufficient condition for Hypothesis (H1) of Theorem 3.1.

The following proposition is a direct consequence of Proposition 5.3 in [18] and provides Razumikhin sufficient conditions that guarantee assumption (H1) for subsystem (3.2). Its proof can be found in the Appendix. Notice that if $V: [-r_1, +\infty) \times \Re^{n_1} \to \Re$ is a locally Lipschitz mapping and $(t, x_1, v) \in \Re^+ \times \Re^{n_1} \times \Re^{n_1}$, we may define the generalized Dini derivative in the direction $v \in \Re^{n_1}$ by the following relation:

$$D^+V(t, x_1; v) := \limsup_{h \to 0^+} \frac{V(t+h, x_1 + hv) - V(t, x_1)}{h} \tag{3.18}$$

**Proposition 3.5 (Razumikhin-like sufficient conditions for hypothesis (H1)):** *Consider system (3.2) under hypotheses (R3-7) and suppose that there exist a locally Lipschitz function $V: [-r_1, +\infty) \times \Re^{n_1} \to \Re^+$, functions $a_1, a_2, a$ of class $K_\infty$ with $a(s) < s$ for all $s > 0$, $\zeta, \zeta^u$ of class $N$, $\beta, \delta_1, \delta_1^u$ of class $K^+$ and a positive definite function $\rho$ such that:*

$$V(t - r_1, x_1) \leq a_2(\beta(t)|x_1|), \quad \forall (t, x_1) \in \Re^+ \times \Re^{n_1} \tag{3.19}$$

$$D^+V(t, x_1(0); f_1(t, d, x_1, u, v_1)) \leq -\rho(V(t, x_1(0))),$$
for all $(t, x_1, u, v_1, d) \in \Re^+ \times C^0([-r_1,0];\Re^{n_1}) \times U \times S_2 \times D$ with

$$\max\left\{ \zeta(\delta_1(t)|v_1|), \zeta^u(\delta_1^u(t)|u|), a\left( \sup_{\theta \in [-r_1,0]} V(t+\theta, x_1(\theta)) \right) \right\} \leq V(t, x_1(0)) \tag{3.20}$$

*Moreover, suppose that there exist functions $a_1, p$ of class $K_\infty$, $\mu$ of class $K^+$ and a constant $R \geq 0$ such that one of the following inequalities holds:*

$$a_1\left(|H(t,x_1)|\right) \leq \sup_{\theta \in [-r_1,0]} V(t+\theta, x_1(\theta)), \quad \forall (t,x_1) \in \Re^+ \times C^0([-r_1,0]; \Re^{n_1}) \tag{3.21a}$$

*or*

$$p(\mu(t)|x_1|) \leq V(t-r_1, x_1) + R, \quad \forall (t,x_1) \in \Re^+ \times \Re^{n_1} \tag{3.21b}$$

*If*

* *(3.21a) holds and there exist functions $\mu_1, c_1, \phi \in K^+$, $g_1, p_1, p_1^u \in N$, such that for every $(t_0, x_{10}, (v_1, u, d)) \in \Re^+ \times C^0([-r_1,0]; \Re^{n_1}) \times L^\infty_{loc}(\Re^+; S_2 \times U \times D)$ the solution $x_1(t)$ of (3.2) with initial condition $T_{r_1}(t_0) x_1 = x_{10}$ corresponding to inputs $(v_1, u, d) \in L^\infty_{loc}(\Re^+; S_2 \times U \times D)$ exists for all $t \geq t_0$ and satisfies (3.16), then there exists a function $\sigma_1 \in KL$, such that estimate (3.4) holds with $\beta_1(t) := \max_{0 \leq \tau \leq t+r_1} \beta(\tau)$, $\gamma_1(s) := a_1^{-1}(\zeta(s))$, $\gamma_1^u(s) := a_1^{-1}(\zeta^u(s))$ for all $(t_0, x_{10}, (v_1, u, d)) \in \Re^+ \times C^0([-r_1,0]; \Re^{n_1}) \times L^\infty_{loc}(\Re^+; S_2 \times U \times D)$ and $t \geq t_0$ for the solution $x_1(t)$ of (3.2) with initial condition $T_{r_1}(t_0) x_1 = x_{10}$ corresponding to inputs $(v_1, u, d) \in L^\infty_{loc}(\Re^+; S_2 \times U \times D)$.*

* *(3.21b) holds and $\delta_1(t) \equiv 1$ then for every $(t_0, x_{10}, (v_1, u, d)) \in \Re^+ \times C^0([-r_1,0]; \Re^{n_1}) \times L^\infty_{loc}(\Re^+; S_2 \times U \times D)$ the solution $x_1(t)$ of (3.2) with initial condition $T_{r_1}(t_0) x_1 = x_{10}$ corresponding to inputs $(v_1, u, d) \in L^\infty_{loc}(\Re^+; S_2 \times U \times D)$ exists for all $t \geq t_0$. Furthermore, for every $\phi \in K^+$ there exist functions $\mu_1, c_1 \in K^+$, $g_1, p_1, p_1^u \in N$, such that (3.16) holds for all $t \geq t_0$. Moreover, if $\phi \in K^+$ is bounded and there exists a constant $L > 0$ such that (3.17) holds then the function $c_1 \in K^+$ is bounded.*

Next sufficient Lyapunov-like conditions for hypothesis (H2) of Theorem 3.1 are presented. The proof of the following theorem can be found at the Appendix.

**Theorem 3.6 (Lyapunov-like sufficient conditions for hypothesis (H2)):** *Consider system (3.3) under hypotheses (R1-5) and suppose that there exist a Lyapunov functional $V: \Re^+ \times L^\infty([-r_2,0]; \Re^{n_2}) \to \Re^+$, functions $a_2$ of class $K_\infty$, $\zeta, \zeta^u$ of class $N$, $\beta, \delta_2, \delta_2^u$ of class $K^+$ and a continuous positive definite function $\rho: \Re^+ \to \Re^+$ such that:*

$$V(t, x_2) \leq a_2\left(\beta(t) \|x_2\|_{r_2}\right), \quad \forall (t, x_2) \in \Re^+ \times L^\infty([-r_2,0]; \Re^{n_2}) \tag{3.22}$$

$$V(t+h, G_h(t, x_2; d, u, v_2)) \leq \max\left\{ \sigma(V(t,x_2), h), \sup_{t \leq \tau \leq t+h} \zeta\left(\delta_2(\tau)|v_2(\tau)|\right), \sup_{t \leq \tau \leq t+h} \zeta^u\left(\delta_2^u(\tau)|u(\tau)|\right) \right\},$$

for all $(t, x_2, u, v_2, d) \in \Re^+ \times L^\infty([-r_2,0]; \Re^{n_2}) \times M_U \times M_{S_1} \times M_D$ and $h \in (0, g(t)]$ \hfill (3.23)

*where*

$$G_h(t, x_2; d, u, v_2) = \begin{cases} x_2(h+\theta), & \theta \in [-r_2, -h] \\ f_2(s, d(s), T_{r_2 - \tau(s)}(-\tau(s))x_2, u(s), v_2(s)); s = t+h+\theta, & \theta \in (-h, 0] \end{cases} \tag{3.24}$$

$$g(t) = \min\left\{1, \min_{t \leq s \leq t+1} \tau(s)\right\} \tag{3.25}$$

*and $\sigma \in KL$ is the function that satisfies*

$$\frac{\partial}{\partial t} \sigma(s,t) = -\rho(\sigma(s,t)) \text{ for all } t, s \geq 0 \tag{3.26a}$$

$$\sigma(s,0) = s \text{ for all } s \geq 0 \tag{3.26b}$$

Moreover, suppose that there exist functions $a_1, p$ of class $K_\infty$, $\mu$ of class $K^+$ and a constant $R \geq 0$ such that one of the following inequalities holds:

$$a_1\left(\left|H_2\left(t, T_{r_2-\tau(t)}(-\tau(t))x_2\right)\right|\right) \leq V(t, x_2), \quad \forall (t, x_2) \in \Re^+ \times L^\infty([-r_2, 0]; \Re^{n_2}) \quad (3.27a)$$

or

$$p(\mu(t)|x_2(0)|) \leq V(t, x_2) + R, \quad \forall (t, x_2) \in \Re^+ \times L^\infty([-r_2, 0]; \Re^{n_2}) \quad (3.27b)$$

If

* (3.27a) holds then there exists a function $\sigma_2 \in KL$, such that estimate (3.6) holds with $\beta_2(t) \equiv \beta(t)$, $\gamma_2(s) := a_1^{-1}(\zeta(s))$, $\gamma_2^u(s) := a_1^{-1}(\zeta^u(s))$ for all $(t_0, x_{20}, (v_2, u, d)) \in \Re^+ \times L^\infty([-r_2, 0]; \Re^{n_2}) \times L^\infty_{loc}(\Re^+; S_1 \times U \times D)$ and $t \geq t_0$ for the solution $x_2(t)$ of (3.3) with initial condition $T_{r_2}(t_0)x_2 = x_{20}$ corresponding to inputs $(v_2, u, d) \in L^\infty_{loc}(\Re^+; S_1 \times U \times D)$.

* (3.27b) holds and $\delta_2(t) \equiv 1$ then for every $\phi \in K^+$ there exist functions $\mu_2, c_2 \in K^+$, $g_2, p_2, p_2^u \in N$, such that the following estimate holds for all $t \geq t_0$

$$\phi(t)\|T_{r_2}(t)x_2\|_{r_2} \leq \max\left\{\mu_2(t-t_0), c_2(t_0), g_2\left(\|x_{20}\|_{r_2}\right), \sup_{t_0 \leq s \leq t} p_2(|v_2(s)|), \sup_{t_0 \leq s \leq t} p_2^u(\delta_2^u(s)|u(s)|)\right\} \quad (3.28)$$

for all $(t_0, x_{20}, (v_2, u, d)) \in \Re^+ \times L^\infty([-r_2, 0]; \Re^{n_2}) \times L^\infty_{loc}(\Re^+; S_1 \times U \times D)$ and $t \geq t_0$ for the solution $x_2(t)$ of (3.3) with initial condition $T_{r_2}(t_0)x_2 = x_{20}$ corresponding to inputs $(v_2, u, d) \in L^\infty_{loc}(\Re^+; S_1 \times U \times D)$. Moreover, if $\phi \in K^+$ is bounded and there exists constant $L > 0$ such that (3.17) holds then the function $c_2 \in K^+$ is bounded.

The following corollary shows how a Lyapunov functional satisfying the assumptions of Theorem 3.7 for system (3.3) can be constructed. Its proof is simple and is omitted.

**Corollary 3.7 (Lyapunov-like sufficient conditions for hypothesis (H2)):** *Consider system (3.3) under hypotheses (R1-5) and suppose that there exist a function $W:[-r_2, +\infty) \times \Re^{n_2} \to \Re^+$, functions $\tilde{a}_1, \tilde{a}_2, b$ of class $K_\infty$, $\zeta, \zeta^u$ of class $N$, $\tilde{\beta}, \delta_2, \delta_2^u$ of class $K^+$ and a constant $\lambda \in [0, 1)$ such that:*

$$W(t-r_2, x_2) \leq \tilde{a}_2\left(\tilde{\beta}(t)|x_2|\right), \quad \forall (t, x_2) \in \Re^+ \times \Re^{n_2} \quad (3.29)$$

$$W(t, f_2(t, d, T_{r_2-\tau(t)}(-\tau(t))x_2, u, v_2)) \leq \max\left\{\lambda \sup_{-r_2 \leq \theta \leq -\tau(t)} W(t+\theta, x_2(\theta)), \zeta(\delta_2(t)|v_2|), \zeta^u(\delta_2^u(t)|u|)\right\}$$

for all $(t, x_2, u, v_2, d) \in \Re^+ \times L^\infty([-r_2, 0]; \Re^{n_2}) \times U \times S_1 \times D$ (3.30)

*Let constant $\mu > 0$ such that $\lambda \exp(\mu r_2) \leq 1$. Define for all $(t, x_2) \in \Re^+ \times L^\infty([-r_2, 0]; \Re^{n_2})$ the functional:*

$$V(t, x_2) = \sup_{-r_2 \leq \theta \leq 0} \exp(\mu\theta) W(t+\theta, x_2(\theta)) \quad (3.31)$$

*Then the functional $V: \Re^+ \times L^\infty([-r_2, 0]; \Re^{n_2}) \to \Re^+$ satisfies inequalities (3.22), (3.23) with $\beta(\tau) := \max_{t \leq s \leq t+r_2} \tilde{\beta}(s)$, $a_2(s) := \tilde{a}_2(s)$, $\sigma(s, t) = s\exp(-\mu t)$ and consequently $\sigma \in KL$ is a function satisfying (3.26a,b) with $\rho(s) := \mu s$. Moreover,*

* *if there exist a function $\tilde{a}_1$ of class $K_\infty$, such that the following inequality holds*

$$\tilde{a}_1\left(\left|H_2(t, T_{r_2-\tau(t)}(-\tau(t))x_2)\right|\right) \leq \sup_{\theta \in [-r_2, 0]} W(t+\theta, x_2(\theta)), \quad \forall (t, x_2) \in \Re^+ \times L^\infty([-r_2, 0]; \Re^{n_2}) \quad (3.32)$$

*then inequality (3.27a) holds with $a_1(s) := \exp(-\mu r_2)\tilde{a}_1(s)$.*

∗ *if there exist functions $p$ of class $K_\infty$, $\mu \in K^+$ and a constant $R \geq 0$ such that the following inequality holds*

$$p(\mu(t)|x_2|) \leq W(t, x_2) + R, \quad \forall (t, x_2) \in \Re^+ \times \Re^{n_2} \tag{3.33}$$

*then inequality (3.27b) holds.*

## 4. Illustrating Examples

The following example shows the applicability of the results of the previous section to nonlinear Neutral Functional Differential Equations.

**Example 4.1:** Consider the following system:

$$\begin{aligned} \dot{x}(t) &= -a\, x(t) + d(t)\varphi(t, \dot{x}(t-2r)) \\ x(t) &\in \Re, \ d(t) \in D = [-1,1] \end{aligned} \tag{4.1}$$

where $a > 0$, $r > 0$ constants and $\varphi : \Re^+ \times \Re \to \Re$ continuous function that satisfies

$$|\varphi(t, x)| \leq c^{-1}|x|, \quad \forall (t, x) \in \Re^+ \times \Re \tag{4.2}$$

for certain constant $c > 1$. Clearly, system (4.1) is a system described by a Neutral Functional Differential Equation of the form (1.12). However, it should be noted that system (4.1) cannot be written in Hale's form (1.10). Consequently, the stability properties of the zero solution of system (4.1) cannot be studied using the results contained in [6]. On the other hand an extension of Theorem 1.6 in [19] can be used with a Lyapunov functional of the form $V(t) = \frac{1}{2}x^2(t) + M \int_{t-2r}^{t} \dot{x}^2(\theta) d\theta$, for appropriate $M > 0$ (since Theorem 1.6 in [19] is not concerned with systems with disturbances). Indeed, by using Theorem 1.6 in [19], it is possible to derive sufficient conditions that guarantee $\lim_{t \to +\infty} |x(t)| = 0$, but further study is required to conclude that $\lim_{t \to +\infty} |\dot{x}(t)| = 0$. Here, we will apply Theorem 3.1 to system (4.1) and we will able to derive sufficient conditions for Robust Global Asymptotic Stability of $x(t)$ as well as of $\dot{x}(t)$.

We have already remarked (in the Introduction) that system (4.1) is associated with the following system described by coupled RFDEs and FDEs:

$$\begin{aligned} \dot{x}_1(t) &= -a\, x_1(t) + d(t)\varphi(t, x_2(t-2r)) \\ x_2(t) &= -a\, x_1(t) + d(t)\varphi(t, x_2(t-2r)) \\ x_1(t) &\in \Re, \ x_2(t) \in \Re, \ d(t) \in D = [-1,1] \end{aligned} \tag{4.3}$$

Specifically, if $T_{2r}(t_0)x_2 = T_{2r}(t_0)\dot{x}$ and $x_1(t_0) = x(t_0)$ then the solution of (4.1) corresponding to input $d \in L^\infty_{loc}(\Re^+; D)$ is related with the solution of (4.3) corresponding to the same input $d \in L^\infty_{loc}(\Re^+; D)$ by the following equations:

$$x_1(t) = x(t), \quad x_2(t) = \dot{x}(t), \text{ for all } t \geq t_0$$

Thus, in order to study the stability properties of the zero solution of system (4.1) we are led to study the stability properties of the zero solution of system (4.3). If we define $r_2 = 2r$, $r_1 = 2r$, $\tau(t) \equiv r$, $H_2(t, x_2) := x_2(-r)$, $H_1(t, x_1) := x_1(0)$, $S_1 = S_2 = \Re$, $H(t, x_1, x_2) := (x_1(0), x_2(-r_2))' \in Y := \Re^2$ for all $(t, x_1, x_2) \in \Re^+ \times C^0([-2r,0];\Re) \times L^\infty([-2r,0];\Re)$, then system (4.3) is of the form (3.1). Furthermore, hypotheses (R1-8) are satisfied. Notice that no external input is present ($U = \{0\}$). Notice that

In order to apply Theorem 3.1, we have to study the stability properties of the zero solution for the two independent subsystems:

$$\dot{x}_1(t) = -a\, x_1(t) + d(t)\varphi(t, v_1(t))$$
$$x_1(t) \in \Re,\ v_1(t) \in \Re,\ d(t) \in D = [-1,1]$$
(4.4)

and

$$\dot{x}_2(t) = -a\, v_2(t) + d(t)\varphi(t, x_2(t - r_2))$$
$$x_2(t) \in \Re,\ v_2(t) \in \Re,\ d(t) \in D = [-1,1]$$
(4.5)

Study of system (4.5):

We define the function

$$W(t, x_2) = |x_2|$$
(4.6)

It should be noted that $W$ satisfies all hypotheses of Corollary 3.7. Particularly, (3.29) holds with $\tilde{a}_2(s) := s$ and $\tilde{\beta}(t) \equiv 1$. Moreover, elementary calculations in conjunction with inequality (4.2) and definition (4.6) give for all $(t, d, v_2, x_2) \in \Re^+ \times [-1,1] \times \Re \times L^\infty([-r_2, 0]; \Re)$ and $\varepsilon_2 \in (0, c-1)$:

$$W(t, -a v_2 + d\varphi(t, x_2(-r_2))) = |-a v_2 + d\varphi(t, x_2(-r_2))| \le a|v_2| + |\varphi(t, x_2(-r_2))|$$
$$\le a|v_2| + c^{-1}|x_2(-r_2)| \le \max\left\{\left(1 + \frac{1}{\varepsilon_2}\right)a|v_2|\ ;\ (1+\varepsilon_2)c^{-1} \sup_{-r_2 \le \theta \le -\tau(t)} W(t+\theta, x_2(\theta))\right\}$$

It follows that $W$ satisfies inequality (3.30) with $\lambda = (1+\varepsilon_2)c^{-1} < 1$, $\zeta(s) := \left(1 + \frac{1}{\varepsilon_2}\right)a\, s$, $\delta_2(t) \equiv 1$, $\zeta^u(s) \equiv 0$ and $\delta_2^u(t) \equiv 1$. Inequalities (3.32), (3.33) hold as well with $\tilde{a}_1(s) = p(s) := s$, $R = 0$ and $\mu(t) \equiv 1$.

By virtue of Theorem 3.6 it follows that there exists a function $\sigma_2 \in KL$, $\mu_2, c_2 \in K^+$, $a_2, p_2, p_2^u \in N$, such that estimates (3.6), (3.7) hold with $\beta_2(t) \equiv 1$, $\gamma_2(s) := c\, a\, \varepsilon_2^{-1}$, $\gamma_2^u(s) := 0$, $\delta_2(t) \equiv 1$, $\delta_2^u(t) = q_2^u(t) \equiv 1$, for all $(t_0, x_{20}, (v_2, u, d)) \in \Re^+ \times L^\infty([-r_2, 0]; \Re^{n_2}) \times L^\infty_{loc}(\Re^+; S_1 \times U \times D)$ and $t \ge t_0$ for the solution $x_2(t)$ of (4.4) with initial condition $T_{r_2}(t_0)x_2 = x_{20}$ corresponding to inputs $(v_2, u, d) \in L^\infty_{loc}(\Re^+; S_1 \times U \times D)$. Moreover, $c_2 \in K^+$ is bounded.

Study of system (4.4):

We define the functional for all $(t, x_1) \in \Re^+ \times C^0([-2r, 0]; \Re)$:

$$V(t, x_1) = \frac{1}{2} x_1^2(0)$$
(4.7)

It should be noted that $V$ satisfies inequalities (3.13), (3.15a,b) with $a_1(s) = a_2(s) = p(s) := \frac{1}{2}s^2$, $R = 0$ and $\beta(t) = \mu(t) \equiv 1$. Moreover, we have for all $(t, x_1, d, v_1) \in \Re^+ \times C^0([-2r, 0]; \Re) \times [-1,1] \times \Re$:

$$V^0(t, x_1; -a x_1(0) + d\varphi(t, v_1)) = -a x_1^2(0) + d\, x_1(0)\varphi(t, v_1)$$

Taking into account (4.2) and completing the squares, we find the following implication which holds for all $(t, x_1, d, v_1) \in \Re^+ \times C^0([-2r, 0]; \Re) \times [-1,1] \times \Re$ and $\varepsilon_1 > 1$:

$$\frac{\varepsilon_1^2}{2a^2c^2}|v_1|^2 \leq V(t,x_1) \Rightarrow V^0(t,x_1;-ax_1(0)+d\varphi(t,v_1)) \leq -2a\left(1-\frac{1}{\varepsilon_1}\right)V(t,x_1)$$

Consequently, inequality (3.14) holds with $\zeta(s) := \frac{\varepsilon_1^2}{2a^2c^2}s^2$, $\rho(s) := 2a\left(1-\frac{1}{\varepsilon_1}\right)s$, $\zeta^u(s) \equiv 0$, $\delta_1(t) \equiv 1$ and $\delta_1^u(t) \equiv 1$. By virtue of Theorem 3.3, it follows that there exist functions $\sigma_1 \in KL$, $\mu_1, c_1 \in K^+$, $a_1, p_1, p_1^u \in N$, such that estimates (3.4), (3.5) hold with $\gamma_1(s) := \frac{\varepsilon_1}{ac}s$, $\gamma_1^u(s) := 0$, $q_1^u(t) \equiv 1$ for all $(t_0, x_{10}, (v_1, u, d)) \in \Re^+ \times C^0([-r_1, 0]; \Re^{n_1}) \times L_{loc}^\infty(\Re^+; S_2 \times U \times D)$ and $t \geq t_0$ for the solution $x_1(t)$ of (4.4) with initial condition $T_{r_1}(t_0)x_1 = x_{10}$ corresponding to inputs $(v_1, u, d) \in L_{loc}^\infty(\Re^+; S_2 \times U \times D)$. Moreover, $c_1 \in K^+$ is bounded.

Study of system (4.3):

The previous analysis shows that Hypotheses (H1), (H2) of Theorem 3.1 are satisfied. Inequality (3.8c) holds with $a(s) = s + \frac{ac}{\varepsilon_1}s$. Using Theorem 3.1 and Remark 3.2(b), we conclude that system (4.3) is Uniformly Robustly Globally Asymptotically Output Stable with output $H(t, x_1, x_2) := (x_1(0), x_2(-r_2))' \in Y := \Re^2$ if the following condition holds:

$$\frac{\varepsilon_1}{\varepsilon_2} < 1$$

Since $\varepsilon_1 > 1$ and $\varepsilon_2 \in (0, c-1)$ are arbitrary, the condition above holds if

$$c > 2 \tag{4.8}$$

Hence, due to the nature of the output map $H(t, x_1, x_2) := (x_1(0), x_2(-r_2))' \in Y := \Re^2$, we are in a position to establish that if (4.8) holds then system (4.3) is Uniformly Robustly Globally Asymptotically Stable. ◁

The following example considers linear time-varying systems described by coupled RFDEs and FDEs.

**Example 4.2:** Consider the linear time-varying system:

$$\begin{aligned}\dot{x}_1(t) &= A(t)x_1(t) + B(t)x_2(t-r) + G_1(t)u(t) \\ x_2(t) &= C(t)x_1(t) + D(t)x_2(t-r) + G_2(t)u(t) \\ x_1(t) &\in \Re^{n_1}, x_2(t) \in \Re^{n_2}, u(t) \in \Re^m, t \geq 0\end{aligned} \tag{4.9}$$

where $r > 0$ and the matrices $A(t), B(t), C(t), D(t), G_1(t), G_2(t)$ have continuous elements. The stability properties of the zero solution for linear systems of the form (4.9) without external inputs (i.e., $u(t) \equiv 0$) have been studied for the autonomous case in [21,24,25]. Here, we study the more general problem of the output stability of system (4.9) with output $Y(t) = x_1(t)$.

Let $c \in K^+$ a non-decreasing function and $\eta \in (0,1)$ such that

$$\max\{|D(t)|, \eta\} \leq c(t), \text{ for all } t \geq 0 \tag{4.10}$$

Define $\phi : [-r, +\infty) \to (0, +\infty)$ by the equation:

$$\phi(t) = \exp\left(-r^{-1}\int_{-r}^{t}\log\left(\eta^{-1}c(s+r)\right)ds\right) \tag{4.11}$$

Notice that $\phi:[-r,+\infty) \to (0,+\infty)$ is non-increasing and that definition (4.11) implies:

$$\phi(t) \leq \frac{\eta}{c(t)}\phi(t-r), \text{ for all } t \geq 0 \tag{4.12}$$

We assume that:

**(S1)** $0 \in \Re^{n_1}$ *is globally asymptotically stable for the system* $\dot{x}_1(t) = A(t)x_1(t)$. *Particularly, there exist a continuously differentiable symmetric positive definite matrix* $P(t) \in \Re^{n_1 \times n_1}$ *and a function* $\mu \in K^+$ *with* $\int_0^{+\infty} \mu(t)dt = +\infty$ *such that:*

$$\dot{P}(t) + P(t)A(t) + A'(t)P(t) \leq -2\mu(t)P(t), \text{ for all } t \geq 0 \tag{4.13a}$$

$$I \leq P(t), \text{ for all } t \geq 0 \tag{4.13b}$$

Notice that hypothesis (S1) **does not** imply that $0 \in \Re^{n_1}$ is **uniformly** globally asymptotically stable for the system $\dot{x}_1(t) = A(t)x_1(t)$. Moreover, if $0 \in \Re^{n_1}$ is globally asymptotically stable for the system $\dot{x}_1(t) = A(t)x_1(t)$ then Proposition 2 in [12] implies the existence of a continuously differentiable symmetric positive definite matrix $P(t) \in \Re^{n_1 \times n_1}$ and a function $\mu \in K^+$ with $\int_0^{+\infty} \mu(t)dt = +\infty$ such that (4.13a,b) hold.

Furthermore, we assume that:

**(S2)** *There exist a constant* $M > 0$ *such that the following inequalities hold:*

$$\sup_{t \geq 0} \frac{|B(t)|\sqrt{|P(t)|}}{\mu(t)\phi(t-r)} \leq M \tag{4.14a}$$

$$\sup_{t \geq 0} \left(\frac{|B(t)|\sqrt{|P(t)|}}{\mu(t)\phi(t-r)} \max_{\tau \in [0,t]} \phi(\tau)|C(\tau)|\right) < 1 - \eta \tag{4.14b}$$

We will next show that system (4.9) with output $Y(t) = x_1(t)$ under hypotheses (S1), (S2) satisfies the WIOS property from the input $u$. First notice that system (4.9) satisfies hypotheses (R1-8) with $r_1 = r_2 = r$, $\tau(t) \equiv r/2$, $H_1(t, x_1) := x_1(0)$, $H_2(t, T_{r-\tau(t)}(-\tau(t))x_2) = \phi(t-r)x_2(t-r)$, $S_1 = Y = \Re^{n_1}$, $S_2 = \Re^{n_2}$ and $U = \Re^m$ ($d \in D$ is irrelevant for this system).

Next consider the system described by linear FDEs:

$$\begin{aligned} x_2(t) &= C(t)v_2(t) + D(t)x_2(t-r) + G_2(t)u(t) \\ v_1(t) &\in \Re^{n_1}, x_2(t) \in \Re^{n_2}, u(t) \in \Re^m, t \geq 0 \end{aligned} \tag{4.15}$$

Notice that the function:

$$W(t, x_2) = \phi(t)|x_2| \tag{4.16}$$

satisfies (3.29) with $\tilde{a}_2(s) := s$, $\tilde{\beta}(t) := \phi(t-r)$. By virtue of (4.10) and (4.12), $W$ satisfies (3.30) with $\lambda = (1+\varepsilon_1)(1+\varepsilon_2)\eta$, $\zeta(s) := (1+\varepsilon_1^{-1})s$, $\zeta^u(s) := (1+\varepsilon_1)(1+\varepsilon_2^{-1})s$, $\delta_2(t) := \phi(t)|C(t)|$, $\delta_2^u(t) := \phi(t)|G_2(t)|$ and $\varepsilon_1, \varepsilon_2 > 0$ with $(1+\varepsilon_1)(1+\varepsilon_2)\eta < 1$. Moreover, inequalities (3.32), (3.33) are satisfied with $\tilde{a}_1(s) := s$, $p(s) := s$, $\mu(t) := \phi(t)$ and $R = 0$. It follows from Theorem 3.6 and Corollary 3.7 that there exist functions $\sigma_2 \in KL$, $\mu_2, c_2 \in K^+$, $a_2, p_2, p_2^u \in N$, such that estimates (3.6), (3.7) hold with $\gamma_2(s) := \eta^{-1}\varepsilon_1^{-1}(1+\varepsilon_2)^{-1}s$, $\gamma_2^u(s) := \eta^{-1}\varepsilon_2^{-1}s$, $\beta_2(t) := \phi(t-r)$, $q_2^u(t) := \phi(t)|G_2(t)|$ for all $(t_0, x_{20}, (v_2, u)) \in \Re^+ \times L^\infty([-r_2, 0]; \Re^{n_2}) \times L^\infty_{loc}(\Re^+; S_1 \times U)$ and $t \geq t_0$ for the solution $x_2(t)$ of (4.15) with initial condition $T_{r_2}(t_0)x_2 = x_{20}$ corresponding to inputs $(v_2, u) \in L^\infty_{loc}(\Re^+; S_2 \times U)$. Moreover, if there exists constant $l > 0$ such that $\phi(t) \geq l$ for all $t \geq 0$ then the function $c_2 \in K^+$ is bounded.

Next, consider the linear system described by RFDEs:

$$\dot{x}_1(t) = A(t)x_1(t) + \frac{1}{\phi(t-r)}B(t)v_1(t) + G_1(t)u(t)$$
$$x_1(t) \in \Re^{n_1}, v_1(t) \in \Re^{n_2}, u(t) \in \Re^m, t \geq 0$$
(4.17)

and define the functional for all $(t, x_1) \in \Re^+ \times C^0([-r, 0]; \Re^{n_1})$:

$$V(t, x_1) := x_1'(0)P(t)x_1(0)$$
(4.18)

It is clear that for every $(t_0, x_{10}) \in \Re^+ \times C^0([-r, 0]; \Re^{n_1})$, $(v_1, u) \in L^\infty_{loc}(\Re^+; \Re^{n_2}) \times L^\infty_{loc}(\Re^+; \Re^m)$, the function $V(t) = V(t, T_r(t)x_1)$ is absolutely continuous on $[t_0, +\infty)$, where $x_1(t)$ denotes the solution of (4.17) with initial condition $T_r(t_0)x_1 = x_{10}$ corresponding to inputs $(v_1, u) \in L^\infty_{loc}(\Re^+; \Re^{n_2}) \times L^\infty_{loc}(\Re^+; \Re^m)$. Using (4.13a) and completing the squares, we can argue that the derivative of $V(t) = V(t, T_r(t)x_1)$ satisfies the following inequalities a.e. on $[t_0, +\infty)$ for all $\varepsilon_3 \in (0,1)$:

$$\dot{V}(t) \leq -2\mu(t)V(t) + \frac{2}{\phi(t-r)}x_1'(t)P(t)B(t)v_1(t) + 2x_1'(t)P(t)G_1(t)u(t)$$

$$\leq -(1-\varepsilon_3)\mu(t)V(t) + \frac{1}{\mu(t)\phi^2(t-r)}v_1'(t)B'(t)P(t)B(t)v_1(t) + \frac{1}{\varepsilon_3\mu(t)}u'(t)G_1'(t)P(t)G_1(t)u(t)$$

The above linear differential inequality implies the following estimate:

$$V(t) \leq V(t_0)\exp\left(-(1-\varepsilon_3)\int_{t_0}^t \mu(s)ds\right) + \int_{t_0}^t \frac{1}{\mu(\tau)\phi^2(\tau-r)}\exp\left(-(1-\varepsilon_3)\int_\tau^t \mu(s)ds\right)|P(\tau)||B(\tau)|^2|v_1(\tau)|^2 d\tau$$

$$+ \int_{t_0}^t \frac{1}{\varepsilon_3\mu(\tau)}\exp\left(-(1-\varepsilon_3)\int_\tau^t \mu(s)ds\right)|P(\tau)||G_1(\tau)|^2|u(\tau)|^2 d\tau$$

$$\leq |P(t_0)|\|x_{10}\|_{r_1}^2 \exp\left((1-\varepsilon_3)\int_0^{t_0} \mu(s)ds\right)\exp\left(-(1-\varepsilon_3)\int_0^{t-t_0} \mu(s)ds\right)$$

$$+ \sup_{t_0 \leq \tau \leq t} \frac{|P(\tau)||B(\tau)|^2}{(1-\varepsilon_3)\mu^2(\tau)\phi^2(\tau-r)}|v_1(\tau)|^2 + \sup_{t_0 \leq \tau \leq t} \frac{|P(\tau)||G_1(\tau)|^2}{\varepsilon_3(1-\varepsilon_3)\mu^2(\tau)}|u(\tau)|^2$$

Consequently, by virtue of (4.13b), estimate (3.4) holds with $\beta_1(t) := \sqrt{|P(t)|} \exp\left(\frac{1-\varepsilon_3}{2} \int_0^t \mu(s)ds\right)$,

$\sigma_1(s,t) := s \exp\left(-\frac{1-\varepsilon_3}{2} \int_0^t \mu(s)ds\right)$, $\gamma_1(s) = \gamma_1^u(s) := s$, $\delta_1(t) := \frac{\sqrt{|P(t)|}\,|B(t)|}{\sqrt{1-\varepsilon_3}\,\mu(t)\phi(t-r)}$ and $\delta_1^u(t) := \frac{\sqrt{|P(t)|}\,|G_1(t)|}{\sqrt{\varepsilon_3(1-\varepsilon_3)}\,\mu(t)}$.

Moreover, there exist functions $\mu_1, c_1 \in K^+$, $g_1, p_1, p_1^u \in N$, such that estimate (3.5) holds with $q_1^u(t) := \frac{\sqrt{|P(t)|}\,|G_1(t)|}{\sqrt{\varepsilon_3(1-\varepsilon_3)}\,\mu(t)}$. Finally, if $\beta_1(t) := \sqrt{|P(t)|} \exp\left(\frac{1-\varepsilon_3}{2} \int_0^t \mu(s)ds\right)$ is bounded then the function $c_1 \in K^+$ is bounded, too.

The previous analysis shows that Hypotheses (H1), (H2) of Theorem 3.1 are satisfied. Inequality (3.8c) holds with $a(s) = s$. By virtue of (4.14a), it follows that (3.8a) holds. Using Theorem 3.1 and Remark 3.2(b), we conclude that system (4.9) satisfies the WIOS property from the input $u$ with output $H(t, x_1, x_2) := x_1(0)$ if the following condition holds:

$$\frac{1}{\varepsilon_1(1+\varepsilon_2)\eta\sqrt{1-\varepsilon_3}} \sup_{t \geq 0}\left(\frac{|B(t)|\sqrt{|P(t)|}}{\mu(t)\phi(t-r)} \max_{\tau \in [0,t]} \phi(\tau)|C(\tau)|\right) < 1$$

Using the facts that $(1+\varepsilon_1)(1+\varepsilon_2)\eta < 1$ and $\varepsilon_3 \in (0,1)$, we conclude that the above condition holds if (4.14b) holds.

Moreover, if there exists constant $L > 0$ such that $\sqrt{|P(t)|} \exp\left(\int_0^t \mu(s)ds\right) + |C(t)| + \frac{1}{\phi(t)} \leq L$, for all $t \geq 0$, then system (4.9) satisfies the UWIOS property from the input $u$ with output $H(t, x_1, x_2) := x_1(0)$.

Consider for example the system:

$$\begin{aligned}
\dot{x}_1(t) &= -\exp(t)x_1(t) + bx_2(t-r) + u(t) \\
x_2(t) &= x_1(t) + 2x_2(t-r) + u(t) \\
x_1(t) &\in \Re, x_2(t) \in \Re, u(t) \in \Re
\end{aligned} \qquad (4.19)$$

where $b \in \Re$ with $|b| < 1$. System (4.19) has the form (4.9) with $A(t) = -\exp(t)$, $C(t) = 1$, $B(t) = b$, $D(t) = 2$, $G_1(t) = G_2(t) = 1$. System (4.19) satisfies hypothesis (S1) with $P(t) = 1$ and $\mu(t) = \exp(t)$. The function $c(t) = 2$ satisfies inequality (4.10) for all $\eta \in (0,1)$. The function $\phi(t)$ defined by (4.11) is given by the equation $\phi(t) = \left(\frac{\eta}{2}\right)^{\frac{t+r}{r}}$. Consequently, hypothesis (S2) will be satisfied if

$\sup_{t \geq 0}\left(\frac{|B(t)|\sqrt{|P(t)|}}{\mu(t)\phi(t-r)}\right) = |b| \sup_{t \geq 0}\left(\exp\left(-t + t r^{-1} \log\left(\frac{2}{\eta}\right)\right)\right) < 1-\eta$. The previous inequality is satisfied for appropriate selection of $\eta \in (0,1)$ if and only if $r > \log(2)$ and $|b| < 1 - 2\exp(-r)$ (delay-dependent condition). Therefore, if $r > \log(2)$ and $|b| < 1 - 2\exp(-r)$ then system (4.19) satisfies the WIOS property from the input $u$ (the weight function is $\delta(t) = K \exp(-st)$ for appropriate constants $K > 0$ and $s \in (0,1)$). Notice that for system (4.19), it can happen that $\limsup_{t \to +\infty} |x_2(t)| = +\infty$. This feature does not disturb our analysis since the WIOS property concerns the output of the system, which is $Y(t) = x_1(t)$.

It should be noted for system (4.9), that if all matrices $A(t), B(t), C(t), D(t), G_1(t), G_2(t)$ are constant, $D$ is Schur stable ($|D| < 1$ is a case frequently studied in the literature), the matrix $P(t) \in \Re^{n_1 \times n_1}$ is constant and the function $\mu(t)$ is constant (i.e., $\mu(t) \equiv \mu > 0$), then hypothesis (S2) guarantees UIOS property from the input $u$. Particularly, in this case hypothesis (S2) takes the form:

$$|C||B|\sqrt{|P|} < \mu(1-|D|)$$

The above condition is in complete agreement with the Linear Matrix Inequalities proposed in [21,28]. ◁

## 5. Conclusions

In this work stability results for systems described by coupled Retarded Functional Differential Equations (RFDEs) and Functional Difference Equations (FDEs) are presented. The motivation for the study of systems described by coupled RFDEs and FDEs is strong, since such systems can be used to study generalized solutions of systems described by neutral functional differential equations and systems described by hyperbolic partial differential equations. The obtained stability results are based on the observation that the composite system can be regarded as the feedback interconnection of a subsystem described by RFDEs and a subsystem described by FDEs. Recent Small-Gain results and Lyapunov-like characterizations of the Weighted Input-to-Output Stability property for systems described by RFDEs and FDEs are employed. Illustrating examples are provided, which show the applicability of the obtained stability results.

## Appendix

**Proof of Fact I in the proof of Theorem 2.3:** Without loss of generality we may assume that the function $a \in K_\infty$ involved in hypothesis (P2) satisfies $a(s) \geq s$ for all $s \geq 0$ and that the function $\beta \in K^+$ involved in hypothesis (P2) is non-decreasing with $\beta(t) \geq 1$ for all $t \geq 0$. Let $\varepsilon > 0$, $T \in \Re^+$ and define:

$$\tilde{\delta}(\varepsilon, T) := \frac{1}{\beta(T+1)} a^{-1}\left( \frac{\exp(-\tilde{L}(\varepsilon, T)-1)}{4\beta(T+1)} a^{-1}\left(\frac{\varepsilon}{9}\right) \right) \tag{A1}$$

where $\tilde{L}(\varepsilon, T)$ is the constant that corresponds to the bounded sets $I := [0, T+1] \subset \Re^+$, $\Omega \subset \{x_1 \in C^0([-r_1, 0]; \Re^{n_1}) : \|x_1\|_{r_1} \leq \varepsilon\} \times \{x_2 \in L^\infty([-r_2, 0]; \Re^{n_2}) : \|x_2\|_{r_2} \leq \varepsilon\} \times \{u \in U : |u| \leq \varepsilon\}$ and satisfies (1.2).

Let $(t_0, x_{10}, x_{20}) \in [0, T] \times C^0([-r_1, 0]; \Re^{n_1}) \times L^\infty([-r_2, 0]; \Re^{n_2})$, $(u, d) \in L^\infty_{loc}(\Re^+; U) \times L^\infty_{loc}(\Re^+; D)$ with $\|x_{10}\|_{r_1} + \|x_{20}\|_{r_2} + \sup_{t \geq 0}|u(t)| < \tilde{\delta}$ (but otherwise arbitrary). By virtue of Theorem 2.1 there exists $t_{\max} \in (t_0, +\infty]$ and a unique pair of mappings $x_1 \in C^0([t_0 - r_1, t_{\max}); \Re^{n_1})$, $x_2 \in L^\infty_{loc}([t_0 - r_2, t_{\max}); \Re^{n_2})$ with $T_{r_1}(t_0)x_1 = x_{10}$, $T_{r_2}(t_0)x_2 = x_{20}$, $x_1$ being absolutely continuous on $[t_0, t_{\max})$ such that (1.1a) holds a.e. for $t \in [t_0, t_{\max})$ and (1.1b) holds for all $t \in (t_0, t_{\max})$. In addition, if $t_{\max} < +\infty$ then for every $M > 0$ there exists $t \in [t_0, t_{\max})$ with $\|T_{r_1}(t)x_1\|_{r_1} > M$.

Define the set:

$$A = \left\{ t \in [t_0, t_{\max}) : \|T_{r_1}(t)x_1\|_{r_1} > \frac{1}{\beta(T+1)} a^{-1}\left(\frac{\varepsilon}{9}\right) \right\} \tag{A2}$$

We distinguish two cases:

1) $A \cap [t_0, t_0 + h] = \varnothing$

2) $A \cap [t_0, t_0 + h] \neq \varnothing$

Case 1:

We will show that (2.1) holds in this case with $\tilde{\delta} := \tilde{\delta}(\varepsilon, T) > 0$ as defined by (A1). If $A \cap [t_0, t_0 + h] = \varnothing$, where $h := h(T) = \min(1; \min\{\tau(s) : s \in [0, T+1]\})$, then $x_1$ is bounded on $[t_0, t_0 + h]$ and consequently we have $t_{\max} > t_0 + h$. Moreover, by virtue of hypothesis (P2) we have for all $t \in [t_0, t_0 + h]$:

$$|x_2(t)| \leq a\left(\beta(t)\|T_{r_1}(t)x_1\|_{r_1}\right) + a\left(\beta(t)\|T_{r_2-\tau(t)}(t-\tau(t))x_2\|_{r_2-\tau(t)}\right) + a(\beta(t)|u(t)|) \tag{A3}$$

Since $A \cap [t_0, t_0 + h] = \varnothing$ (which implies $a\left(\beta(t)\|T_{r_1}(t)x_1\|_{r_1}\right) \leq \frac{\varepsilon}{9}$ for all $t \in [t_0, t_0 + h]$; see (A2)), $\|x_{20}\|_{r_2} < \tilde{\delta}$ (which implies $a\left(\beta(t)\|T_{r_2-\tau(t)}(t-\tau(t))x_2\|_{r_2-\tau(t)}\right) \leq \frac{\varepsilon}{9}$ for all $t \in [t_0, t_0 + h]$; (see (A1)) and $\sup_{t \geq 0}|u(t)| < \tilde{\delta}$ (which implies $a(\beta(t)|u(t)|) \leq \frac{\varepsilon}{9}$ for all $t \in [t_0, t_0 + h]$; see (A1)), we obtain from (A3):

$$|x_2(t)| \leq \frac{\varepsilon}{3}, \quad \forall t \in [t_0, t_0 + h] \tag{A4}$$

Inequality (A4) in conjunction with the fact that $\|x_{20}\|_{r_2} < \tilde{\delta} \leq \frac{\varepsilon}{3}$ and the assumption $A \cap [t_0, t_0 + h] = \varnothing$ (which implies $\|T_{r_1}(t)x_1\|_{r_1} \leq a\left(\beta(t)\|T_{r_1}(t)x_1\|_{r_1}\right) \leq \frac{\varepsilon}{9}$ for all $t \in [t_0, t_0 + h]$) shows that (2.1) holds in this case.

Case 2:

We will show that this case cannot happen by contradiction. Assume that $A \cap [t_0, t_0 + h] \neq \varnothing$ and define $t_1 = \inf A$. By continuity of the map $t \to \|T_{r_1}(t)x_1\|_{r_1}$ and since $\|x_{10}\|_{r_1} < \tilde{\delta} < \frac{1}{\beta(T+1)} a^{-1}\left(\frac{\varepsilon}{9}\right)$, it follows that $t_1 > t_0$. Hence, by continuity of the map $t \to \|T_{r_1}(t)x_1\|_{r_1}$ and definition (A2) we have $\|T_{r_1}(t_1)x_1\|_{r_1} = \frac{1}{\beta(T+1)} a^{-1}\left(\frac{\varepsilon}{9}\right)$. Evaluating the derivative of the absolutely continuous map $|x_1(t)|^2$ on $[t_0, t_1]$ in conjunction with hypothesis (P4) gives:

$$\frac{d}{dt}|x_1(t)|^2 = 2x_1'(t) f_1(t, d(t), T_{r_1}(t)x_1, T_{r_2-\tau(t)}(t-\tau(t))x_2, u(t))$$

$$\leq 2\tilde{L}\|T_{r_1}(t)x_1\|_{r_1}^2 + 2x_1'(t) f_1(t, d(t), 0, T_{r_2-\tau(t)}(t-\tau(t))x_2, u(t)) \tag{A5}$$

where $\tilde{L}$ is the constant that corresponds to the bounded sets $I := [0, T+1] \subset \Re^+$, $\Omega \subset \left\{ x_1 \in C^0([-r_1, 0]; \Re^{n_1}) : \|x_1\|_{r_1} \leq \varepsilon \right\} \times \left\{ x_2 \in L^\infty([-r_2, 0]; \Re^{n_2}) : \|x_2\|_{r_2} \leq \varepsilon \right\} \times \left\{ u \in U : |u| \leq \varepsilon \right\}$ and satisfies (1.2). Inequality (A5) in conjunction with hypothesis (P2) and the facts that $\|x_{20}\|_{r_2} < \tilde{\delta}$ (which implies

$$a\!\left( \beta(t) \|T_{r_2 - \tau(t)}(t - \tau(t)) x_2 \|_{r_2 - \tau(t)} \right) \leq \frac{\exp(-\tilde{L}(\varepsilon, T) - 1)}{4\beta(T+1)} a^{-1}\!\left(\frac{\varepsilon}{9}\right)$$

for all $t \in [t_0, t_0 + h]$; (see (A1)) and $\sup_{t \geq 0} |u(t)| < \tilde{\delta}$ (which implies $a(\beta(t)|u(t)|) \leq \frac{\exp(-\tilde{L}(\varepsilon, T) - 1)}{4\beta(T+1)} a^{-1}\!\left(\frac{\varepsilon}{9}\right)$ for all $t \in [t_0, t_0 + h]$; see (A1)), gives:

$$\frac{d}{dt}|x_1(t)|^2 \leq 2(\tilde{L}+1)\|T_{r_1}(t) x_1\|_{r_1}^2 + \frac{\exp(-2\tilde{L}-2)}{4\beta^2(T+1)}\left( a^{-1}\!\left(\frac{\varepsilon}{9}\right) \right)^2 \tag{A6}$$

Integrating both sides of (A6) we get for all $t \in [t_0, t_1]$:

$$|x_1(t)|^2 \leq |x_1(t_0)|^2 + 2(\tilde{L}+1)\int_{t_0}^{t} \|T_{r_1}(s) x_1\|_{r_1}^2 ds + \frac{\exp(-2\tilde{L}-2)}{4\beta^2(T+1)}\left( a^{-1}\!\left(\frac{\varepsilon}{9}\right) \right)^2 \tag{A7}$$

The following inequality is a direct consequence of (A7) and holds for all $t \in [t_0, t_1]$:

$$\|T_{r_1}(t) x_1\|_{r_1}^2 \leq \|T_{r_1}(t_0) x_1\|_{r_1}^2 + 2(\tilde{L}+1)\int_{t_0}^{t} \|T_{r_1}(s) x_1\|_{r_1}^2 ds + \frac{\exp(-2\tilde{L}-2)}{4\beta^2(T+1)}\left( a^{-1}\!\left(\frac{\varepsilon}{9}\right) \right)^2 \tag{A8}$$

Since the map $t \to \|T_{r_1}(t) x_1\|_{r_1}$ is continuous and (A8) holds on $[t_0, t_1]$, we may apply the Gronwall-Bellman lemma. We obtain for all $t \in [t_0, t_1]$:

$$\|T_{r_1}(t) x_1\|_{r_1} \leq \exp\!\left( (\tilde{L}+1) \right) \|T_{r_1}(t_0) x_1\|_{r_1} + \frac{1}{2\beta(T+1)} a^{-1}\!\left(\frac{\varepsilon}{9}\right) \tag{A9}$$

Since $\|x_{10}\|_{r_1} < \tilde{\delta}$ (which implies $\|T_{r_1}(t_0) x_1\|_{r_1} \leq \frac{\exp(-\tilde{L}-1)}{4\beta(T+1)} a^{-1}\!\left(\frac{\varepsilon}{9}\right)$) we get from (A9):

$$\|T_{r_1}(t_1) x_1\|_{r_1} \leq \frac{3}{4\beta(T+1)} a^{-1}\!\left(\frac{\varepsilon}{9}\right) < \frac{1}{\beta(T+1)} a^{-1}\!\left(\frac{\varepsilon}{9}\right) \tag{A10}$$

Inequality (A10) contradicts the equality $\|T_{r_1}(t_1) x_1\|_{r_1} = \frac{1}{\beta(T+1)} a^{-1}\!\left(\frac{\varepsilon}{9}\right)$. Thus the case $A \cap [t_0, t_0 + h] \neq \varnothing$ cannot happen.

The proof is complete. ◁

**Definition of the notion of Weighted Input-to-Output Stability (see also [15,16,17]):** *Consider a control system $\Sigma := (X, Y, M_U, M_D, \phi, \pi, H)$ with outputs and the BIC property and for which $0 \in X$ is a robust equilibrium point from the input $u \in M_U$ (in the sense described in [15]). Suppose that $\Sigma$ is Robustly Forward Complete from the input $u \in M_U$ (in the sense described in [15]).*

∗ *If there exist functions $\sigma \in KL$, $\beta, \delta \in K^+$, $\gamma \in N$ such that the following estimate holds for all $u \in M_U$, $(t_0, x_0, d) \in \Re^+ \times X \times M_D$ and $t \geq t_0$:*

$$\|H(t,\phi(t,t_0,x_0,u,d),u(t))\|_Y \leq \sigma\big(\beta(t_0)\|x_0\|_X, t-t_0\big) + \sup_{t_0 \leq \tau \leq t} \gamma\big(\delta(\tau)\|u(\tau)\|_U\big) \quad (A11)$$

*then we say that* $\Sigma$ *satisfies the Weighted Input-to-Output Stability (WIOS) property from the input* $u \in M_U$ *with gain* $\gamma \in N$ *and weight* $\delta \in K^+$. *Moreover, if* $\beta(t) \equiv 1$ *then we say that* $\Sigma$ *satisfies the Uniform Weighted Input-to-Output Stability (UWIOS) property from the input* $u \in M_U$ *with gain* $\gamma \in N$ *and weight* $\delta \in K^+$.

\* *If there exist functions* $\sigma \in KL$, $\beta \in K^+$, $\gamma \in N$ *such that the following estimate holds for all* $u \in M_U$, $(t_0, x_0, d) \in \Re^+ \times X \times M_D$ *and* $t \geq t_0$:

$$\|H(t,\phi(t,t_0,x_0,u,d),u(t))\|_Y \leq \sigma\big(\beta(t_0)\|x_0\|_X, t-t_0\big) + \sup_{t_0 \leq \tau \leq t} \gamma\big(\|u(\tau)\|_U\big) \quad (A12)$$

*then we say that* $\Sigma$ *satisfies the Input-to-Output Stability (IOS) property from the input* $u \in M_U$ *with gain* $\gamma \in N$. *Moreover, if* $\beta(t) \equiv 1$ *then we say that* $\Sigma$ *satisfies the Uniform Input-to-Output Stability (UIOS) property from the input* $u \in M_U$ *with gain* $\gamma \in N$.

*For the special case of the identity output mapping, i.e.,* $H(t,x,u) := x$, *the (Uniform) (Weighted) Input-to-Output Stability property from the input* $u \in M_U$ *is called (Uniform) (Weighted) Input-to-State Stability property from the input* $u \in M_U$.

*Finally, if no external input is present* ($U = \{0\}$) *then we say that* $\Sigma$ *is (Uniformly) Robustly Globally Asymptotically Output Stable ((U)RGAOS). For the special case of the identity output mapping, i.e.,* $H(t,x,u) := x$, *we say that* $\Sigma$ *is (Uniformly) Robustly Globally Asymptotically Stable ((U)RGAS).*

**Proof of Theorem 3.3:** Let arbitrary $(t_0, x_{10}, (v_1, u, d)) \in \Re^+ \times C^0([-r_1, 0]; \Re^{n_1}) \times L^\infty_{loc}(\Re^+; S_2 \times U \times D)$ and consider the solution $x_1(t)$ of (3.2) with initial condition $T_{r_1}(t_0) = x_{10}$ corresponding to inputs $(v_1, u, d) \in L^\infty_{loc}(\Re^+; S_2 \times U \times D)$. By virtue of Theorem 4.6 in [18], implication (3.14) guarantees the existence of a function $\sigma \in KL$ with $\sigma(s,0) = s$ for all $s \geq 0$, such that the following estimate holds for all $t \geq t_0$:

$$V(t, T_{r_1}(t)x_1) \leq \max\left\{ \sigma\big(a_2(\beta(t_0)\|x_{10}\|_{r_1}), t-t_0\big), \sup_{t_0 \leq \tau \leq t} \sigma\big(\zeta(\delta_1(\tau)|v_1(\tau)|), t-\tau\big), \sup_{t_0 \leq \tau \leq t} \sigma\big(\zeta^u(\delta_1^u(\tau)|u(\tau)|), t-\tau\big) \right\}$$
(A13)

We next distinguish the following cases:

**A)** If (3.15a) holds then, by exploiting the left hand-side inequality (3.14) and the fact $\sigma \in KL$ with $\sigma(s,0) = s$ for all $s \geq 0$, we obtain from (A13) the following estimate, which holds for all $t \geq t_0$:

$$a_1\big(|H_1(t,T_{r_1}(t)x_1)|\big) \leq \max\left\{ \sigma\big(a_2(\beta(t_0)\|x_{10}\|_{r_1}), t-t_0\big), \sup_{t_0 \leq \tau \leq t} \zeta(\delta_1(\tau)|v_1(\tau)|), \sup_{t_0 \leq \tau \leq t} \zeta^u(\delta_1^u(\tau)|u(\tau)|) \right\} \quad (A14)$$

Estimate (A14) implies that estimate (3.4) holds with $\beta_1(t) \equiv \beta(t)$, $\gamma_1(s) := a_1^{-1}(\zeta(s))$, $\gamma_1^u(s) := a_1^{-1}(\zeta^u(s))$.

**B)** If (3.15b) holds and $\delta_1(t) \equiv 1$ then Theorem 4.6 in [18] directly implies that system (3.2) is Robustly Forward Complete (see [15]) and consequently for all $(t_0, x_{10}, (v_1, u, d)) \in \Re^+ \times C^0([-r_1, 0]; \Re^{n_1}) \times L^\infty_{loc}(\Re^+; S_2 \times U \times D)$ the solution $x_1(t)$ of (3.2) with initial condition $T_{r_1}(t_0)x_1 = x_{10}$ corresponding to inputs $(v_1, u, d) \in L^\infty_{loc}(\Re^+; S_2 \times U \times D)$ exists for all $t \geq t_0$. By exploiting the left hand-side inequality (3.15b) and the fact $\sigma \in KL$ with $\sigma(s,0) = s$ for all $s \geq 0$, we obtain from (A13) the following estimate, which holds for all $t \geq t_0$:

$$\mu(t)|x_1(t)| \leq p^{-1}(2R) + \max\left\{ p^{-1}\left(2a\left(\beta(t_0)\|x_{10}\|_{r_1}\right)\right), \sup_{t_0 \leq \tau \leq t} p^{-1}\left(2\zeta(|v_1(\tau)|)\right), \sup_{t_0 \leq \tau \leq t} p^{-1}\left(2\zeta^u\left(\delta_1^u(\tau)|u(\tau)|\right)\right) \right\} \quad (A15)$$

By virtue of Corollary 10 and Remark 11 in [33] there exists $\kappa \in K_\infty$ such that $p^{-1}(2a(rs)) \leq \kappa(r)\kappa(s)$ for all $(r,s) \in (\Re^+)^2$. Consequently, we obtain from (A15):

$$\mu(t)|x_1(t)| \leq p^{-1}(2R) + \frac{1}{2}(\kappa(\beta(t_0)))^2 + \max\left\{ \frac{1}{2}\left(\kappa(\|x_{10}\|_{r_1})\right)^2, \sup_{t_0 \leq \tau \leq t} p^{-1}\left(2\zeta(|v_1(\tau)|)\right), \sup_{t_0 \leq \tau \leq t} p^{-1}\left(2\zeta^u\left(\delta_1^u(\tau)|u(\tau)|\right)\right) \right\}$$

which implies:

$$|x_1(t)| \leq \frac{1}{2\mu^2(t)} + \left( p^{-1}(2R) + \frac{1}{2}(\kappa(\beta(t_0)))^2 \right)^2 \\ + \max\left\{ \left(\kappa(\|x_{10}\|_{r_1})\right)^4, \sup_{t_0 \leq \tau \leq t} \left(p^{-1}\left(2\zeta(|v_1(\tau)|)\right)\right)^2, \sup_{t_0 \leq \tau \leq t} \left(p^{-1}\left(2\zeta^u\left(\delta_1^u(\tau)|u(\tau)|\right)\right)\right)^2 \right\} \quad (A16)$$

Estimate (A16) shows that the following estimate holds for all $t \geq t_0$:

$$\|T_{r_1}(t)x_1\|_{r_1} \leq \frac{1}{2\min\{\mu^2(\tau): \tau \in [0,t]\}} + \left( p^{-1}(2R) + \frac{1}{2}(\kappa(\beta(t_0)))^2 \right)^2 \\ + \max\left\{ \left(\kappa(\|x_{10}\|_{r_1})\right)^4 + \|x_{10}\|_{r_1}, \sup_{t_0 \leq \tau \leq t} \left(p^{-1}\left(2\zeta(|v_1(\tau)|)\right)\right)^2, \sup_{t_0 \leq \tau \leq t} \left(p^{-1}\left(2\zeta^u\left(\delta_1^u(\tau)|u(\tau)|\right)\right)\right)^2 \right\}$$

Let $\phi \in K^+$. Multiplying the above inequality by $\phi(t)$ and using repeatedly the inequality $ab \leq \frac{1}{2}a^2 + \frac{1}{2}b^2$ we get:

$$\phi(t)\|T_{r_1}(t)x_1\|_{r_1} \leq \frac{\phi(t)}{2\min\{\mu^2(\tau):\tau\in[0,t]\}} + \phi^2(t) + \frac{1}{2}\left( p^{-1}(2R) + \frac{1}{2}(\kappa(\beta(t_0)))^2 \right)^4 \\ + \max\left\{ \left(\left(\kappa(\|x_{10}\|_{r_1})\right)^4 + \|x_{10}\|_{r_1}\right)^2, \sup_{t_0 \leq \tau \leq t} \left(p^{-1}\left(2\zeta(|v_1(\tau)|)\right)\right)^4, \sup_{t_0 \leq \tau \leq t} \left(p^{-1}\left(2\zeta^u\left(\delta_1^u(\tau)|u(\tau)|\right)\right)\right)^4 \right\}$$

Using the fact that $a + b \leq \max\{2a, 2b\}$ for all $a, b \geq 0$ in conjunction with the above inequality gives:

$$\phi(t)\|T_{r_1}(t)x_1\|_{r_1} \leq \max\left\{ \begin{array}{l} \dfrac{\phi(t)}{\min\{\mu^2(\tau):\tau\in[0,t]\}} + 2\phi^2(t) + \left( p^{-1}(2R) + \frac{1}{2}(\kappa(\beta(t_0)))^2 \right)^4, 2\left(\left(\kappa(\|x_{10}\|_{r_1})\right)^4 + \|x_{10}\|_{r_1}\right)^2 \\ \sup_{t_0 \leq \tau \leq t} 2\left(p^{-1}\left(2\zeta(|v_1(\tau)|)\right)\right)^4, \sup_{t_0 \leq \tau \leq t} 2\left(p^{-1}\left(2\zeta^u\left(\delta_1^u(\tau)|u(\tau)|\right)\right)\right)^4 \end{array} \right\} \quad (A17)$$

Notice that if (3.17) holds and $\phi \in K^+$ is bounded then estimate (A17) shows that (3.16) holds for appropriate $\mu_1, c_1 \in K^+$, $g_1, p_1, p_1^u \in \mathcal{N}$ with $c_1 \in K^+$ being bounded.

Let $\gamma(t) = \dfrac{\phi(t)}{\min\{\mu^2(\tau):\tau\in[0,t]\}} + 2\phi^2(t)$. Define:

$$a(T,s) := \max\{\gamma(t_0 + h) - \gamma(t_0) : h \in [0,s], t_0 \in [0,T]\} \quad (A18)$$

Clearly, definition (A18) implies that for each fixed $s \geq 0$ $a(\cdot, s)$ is non-decreasing and for each fixed $T \geq 0$ $a(T, \cdot)$ is non-decreasing. Furthermore, continuity of $\gamma$ guarantees that for every $T \geq 0$ $\lim_{s \to 0^+} a(T,s) = a(T,0) = 0$. It turns out from Lemma 2.3 in [11], that there exist functions $\zeta \in K_\infty$ and $q \in K^+$ such that

$$a(T,s) \leq \zeta(q(T)s), \quad \forall (T,s) \in (\Re^+)^2 \tag{A19}$$

Combining definition (A18) with inequality (A19), we conclude that for all $t_0 \geq 0$ and $t \geq t_0$, it holds that:

$$\gamma(t) \leq \gamma(t_0) + \zeta\big(q(t_0)(t-t_0)\big) \leq \gamma(t_0) + \zeta\left(\frac{1}{2}q^2(t_0) + \frac{1}{2}(t-t_0)^2\right)$$

$$\leq \gamma(t_0) + \zeta\big(q^2(t_0)\big) + \zeta\big((t-t_0)^2\big) \leq \max\left\{2\gamma(t_0) + 2\zeta\big(q^2(t_0)\big); 2\zeta\big((t-t_0)^2\big)\right\}$$

The above inequality in conjunction with (A17) and definition $\gamma(t) = \dfrac{\phi(t)}{\min\{\mu^2(\tau) : \tau \in [0,t]\}} + 2\phi^2(t)$ implies that (3.16) holds for appropriate $\mu_1, c_1 \in K^+$, $g_1, p_1, p_1^u \in N$. The proof is complete. ◁

**Proof of Proposition 3.5:** Let arbitrary $(t_0, x_{10}, (v_1, u, d)) \in \Re^+ \times C^0([-r_1, 0]; \Re^{n_1}) \times L_{loc}^\infty(\Re^+; S_2 \times U \times D)$ and consider the solution $x_1(t)$ of (3.2) with initial condition $T_{r_1}(t_0)x_1 = x_{10}$ corresponding to inputs $(v_1, u, d) \in L_{loc}^\infty(\Re^+; S_2 \times U \times D)$. Let $t_{\max} > t_0$ be the maximal existence time of the solution. It follows from (3.20) and Lemma 4.7 in [18] that there exists a continuous function $\sigma$ of class $KL$, with $\sigma(s, 0) = s$ for all $s \geq 0$ such that for all $t \in [t_0, t_{\max})$ we have:

$$V(t, x_1(t)) \leq \max\left\{\begin{array}{l} \sigma\big(V(t_0, x_1(t_0)), t-t_0\big); \sup_{t_0 \leq s \leq t} \sigma\big(\zeta(\delta_1(s)|v_1(s)|), t-s\big) \\ \sup_{t_0 \leq s \leq t} \sigma\left(a\left(\sup_{\theta \in [-r_1, 0]} V(s+\theta, x_1(s+\theta))\right), t-s\right); \sup_{t_0 \leq s \leq t} \sigma\big(\zeta^u(\delta_1^u(s)|u(s)|), t-s\big) \end{array}\right\} \tag{A20}$$

An immediate consequence of estimate (A20) and the fact that $\sigma(s,0) = s$ for all $s \geq 0$ is the following estimate for all $t \in [t_0, t_{\max})$:

$$\sup_{\theta \in [-r_1, 0]} V(t+\theta, x_1(t+\theta)) \leq \max\left\{\begin{array}{l} \overline{\sigma}\left(\sup_{\theta \in [-r_1, 0]} V(t_0+\theta, x_1(t_0+\theta)), t-t_0\right); \sup_{t_0 \leq s \leq t} \zeta\big(\delta_1(s)|v_1(s)|\big) \\ a\left(\sup_{t_0 \leq s \leq t} \sup_{\theta \in [-r_1, 0]} V(s+\theta, x_1(s+\theta))\right); \sup_{t_0 \leq s \leq t} \zeta^u\big(\delta_1^u(s)|u(s)|\big) \end{array}\right\} \tag{A21}$$

where $\overline{\sigma}(s,t) := s$ for $t \in [0, r]$ and $\overline{\sigma}(s,t) := \sigma(s, t-r)$ for $t > r$. Using the fact that $a(s) < s$ for all $s > 0$ and estimate (A21) it may be shown that:

$$\sup_{\theta \in [-r_1, 0]} V(t+\theta, x_1(t+\theta)) \leq \max\left\{\sup_{\theta \in [-r_1, 0]} V(t_0+\theta, x_1(t_0+\theta)); \sup_{t_0 \leq s \leq t} \zeta\big(\delta_1(s)|v_1(s)|\big); \sup_{t_0 \leq s \leq t} \zeta^u\big(\delta_1^u(s)|u(s)|\big)\right\}, \quad t \in [t_0, t_{\max}) \tag{A22}$$

Combining (A21) with (A22) we obtain for all $t \in [t_0, t_{\max})$:

$$\sup_{\theta \in [-r_1, 0]} V(t+\theta, x_1(t+\theta)) \leq \inf_{t_0 \leq \xi \leq t} \max\left\{\begin{array}{l} \overline{\sigma}\left(\sup_{\theta \in [-r_1, 0]} V(t_0+\theta, x_1(t_0+\theta)), t-\xi\right); a\left(\sup_{\xi \leq s \leq t} \sup_{\theta \in [-r_1, 0]} V(s+\theta, x_1(s+\theta))\right) \\ \sup_{t_0 \leq s \leq t} \zeta\big(\delta_1(s)|v_1(s)|\big); \sup_{t_0 \leq s \leq t} \zeta^u\big(\delta_1^u(s)|u(s)|\big) \end{array}\right\} \tag{A23}$$

Theorem 1 in [35] in conjunction with inequality (A23) implies the existence of $\tilde{\sigma} \in KL$ such that:

$$\sup_{\theta \in [-r_1, 0]} V(t+\theta, x_1(t+\theta)) \leq \max \left\{ \begin{array}{l} \tilde{\sigma}\left( \sup_{\theta \in [-r_1, 0]} V(t_0+\theta, x_1(t_0+\theta)), t-t_0 \right) \\ \sup_{t_0 \leq s \leq t} \zeta(\delta_1(s)|v_1(s)|); \sup_{t_0 \leq s \leq t} \zeta^u(\delta_1^u(s)|u(s)|) \end{array} \right\}, \quad \forall t \in [t_0, t_{\max}) \quad (A24)$$

We next distinguish the following cases:

**A)** If (3.21a) holds and there exist functions $\mu_1, c_1, \phi \in K^+$, $g_1, p_1, p_1^u \in N$, such that for every $(t_0, x_{10}, (v_1, u, d)) \in \Re^+ \times C^0([-r_1, 0]; \Re^{n_1}) \times L^\infty_{loc}(\Re^+; S_2 \times U \times D)$ the solution $x_1(t)$ of (3.2) with initial condition $T_{r_1}(t_0)x_1 = x_{10}$ corresponding to inputs $(v_1, u, d) \in L^\infty_{loc}(\Re^+; S_2 \times U \times D)$ exists for all $t \geq t_0$ and satisfies (3.16), then we clearly have $t_{\max} = +\infty$. Inequalities (3.19), (3.21a) in conjunction with estimate (A24), guarantee that there exist a function $\sigma_1 \in KL$, such that estimate (3.4) holds with $\gamma_1(s) := a_1^{-1}(\zeta(s))$, $\gamma_1^u(s) := a_1^{-1}(\zeta^u(s))$, $\beta_1(t) := \max_{0 \leq \tau \leq t+r_1} \beta(\tau)$, for all $(t_0, x_{10}, (v_1, u, d)) \in \Re^+ \times C^0([-r_1, 0]; \Re^{n_1}) \times L^\infty_{loc}(\Re^+; S_2 \times U \times D)$ and $t \geq t_0$ for the solution $x_1(t)$ of (3.2) with initial condition $T_{r_1}(t_0) = x_{10}$ corresponding to inputs $(v_1, u, d) \in L^\infty_{loc}(\Re^+; S_2 \times U \times D)$.

**B)** If (3.21b) holds and $\delta_1(t) \equiv 1$, then using (3.19) and (3.21b), we get for all $t \in [t_0, t_{\max})$:

$$\|T_{r_1}(t)x_1\|_{r_1} \leq \frac{1}{2 \min_{0 \leq \tau \leq t+r} \mu^2(\tau)} + \left(p^{-1}(2R)\right)^2$$

$$+ \max \left\{ \left(p^{-1}\left(2a_2\left(\max_{0 \leq \tau \leq t_0+r_1} \beta(\tau) \|x_{10}\|_{r_1}\right)\right)\right)^2 ; \sup_{t_0 \leq s \leq t} \left(p^{-1}(2\zeta(|v_1(s)|))\right)^2 ; \sup_{t_0 \leq s \leq t} \left(p^{-1}(2\zeta^u(\delta_1^u(s)|u(s)|))\right)^2 \right\} \quad (A25)$$

It follows from estimate (A25) and a simple contradiction argument that for all $(t_0, x_{10}, (v_1, u, d)) \in \Re^+ \times C^0([-r_1, 0]; \Re^{n_1}) \times L^\infty_{loc}(\Re^+; S_2 \times U \times D)$ the solution $x_1(t)$ of (3.2) with initial condition $T_{r_1}(t_0)x_1 = x_{10}$ corresponding to inputs $(v_1, u, d) \in L^\infty_{loc}(\Re^+; S_2 \times U \times D)$ exists for all $t \geq t_0$, i.e., $t_{\max} = +\infty$. From this point the proof continues in exactly the same way as in Case B in the proof of Theorem 3.3. ◁

**Proof of Theorem 3.6:** By virtue of Lemma 4.4 in [20] there exists $\sigma \in KL$ satisfying (3.26a,b). Let arbitrary $(t_0, x_{20}, (v_2, u, d)) \in \Re^+ \times L^\infty([-r_2, 0]; \Re^{n_2}) \times L^\infty_{loc}(\Re^+; S_1 \times U \times D)$ and consider the solution $x_2(t)$ of (3.3) with initial condition $T_{r_2}(t_0)x_2 = x_{20}$ corresponding to inputs $(v_2, u, d) \in L^\infty_{loc}(\Re^+; S_1 \times U \times D)$. Define the sequence $\{t_i\}_{i=0}^\infty$:

$$t_{i+1} = t_i + g(t_i) \quad (A26)$$

Working exactly as in the proof of Theorem 2.1, it can be shown by contradiction that $\lim t_i = +\infty$. Moreover, using induction, (3.26b), inequality (3.23), as well as the fact that $\sigma(\sigma(s,t),h) = \sigma(s, t+h)$ for all $s, t, h \geq 0$, we get for all $i = 0, 1, 2, \ldots$:

$$V(t_i, T_{r_2}(t_i)x_2) \leq \max \left\{ \sigma(V(t_0, x_{20}), t_i - t_0), \sup_{t_0 \leq \tau \leq t_i} \zeta(\delta_2(\tau)|v_2(\tau)|), \sup_{t_0 \leq \tau \leq t_i} \zeta^u(\delta_2^u(\tau)|u(\tau)|) \right\} \quad (A27)$$

Using (A32), (3.26b), inequality (3.23), as well as the fact that $\sigma(\sigma(s,t),h) = \sigma(s, t+h)$ for all $s, t, h \geq 0$, we can establish that the following estimate holds for all $t \notin \{t_i\}_{i=0}^\infty$:

$$V(t, T_{r_2}(t)x_2) \le \max\left\{\sigma(V(t_0, x_{20}), t-t_0), \sup_{t_0 \le \tau \le t} \zeta(\delta_2(\tau)|v_2(\tau)|), \sup_{t_0 \le \tau \le t} \zeta^u(\delta_2^u(\tau)|u(\tau)|)\right\} \quad \text{(A28)}$$

By virtue of (A27), (A28), we conclude that (A28) holds for all $t \ge t_0$. Next we distinguish the cases:

**A)** If (3.27a) holds, then by combining (A28) with (3.22) and (3.27a), we conclude that estimate (3.6) holds with $\gamma_2(s) := a_1^{-1}(\zeta(s))$, $\gamma_2^u(s) := a_1^{-1}(\zeta^u(s))$ and $\sigma_2(s,t) = a_1^{-1}(\sigma(a_2(s), t))$.

**B)** If (3.27b) holds and $\delta_2(t) \equiv 1$, then estimate (A28) in conjunction with (3.22) and (3.27b) implies the following estimate:

$$p(\mu(t)|x_2(t)|) \le R + \max\left\{a_2(\beta(t_0)\|x_{20}\|_{r_2}), \sup_{t_0 \le \tau \le t} \zeta(|v_2(\tau)|), \sup_{t_0 \le \tau \le t} \zeta^u(\delta_2^u(\tau)|u(\tau)|)\right\}$$

From this point the proof continues in exactly the same way as in Case B in the proof of Theorem 3.3. ◁